\documentclass[final,onefignum,onetabnum]{siamart190516}

\clearpage{}

\usepackage{bm}

\newcommand{\RR}{\mathbb{R}}
\newcommand{\CC}{\mathbb{C}}

\newcommand{\eps}{\varepsilon}

\newcommand{\set}[1]{\mathsf{#1}}
\newcommand{\vct}[1]{\bm{#1}}
\newcommand{\mtx}[1]{\bm{#1}}

\newcommand{\Id}{\mathbf{I}}

\newcommand{\lspan}{\operatorname{span}}

\newcommand{\range}{\operatorname{range}}

\newcommand{\abs}[1]{\vert #1 \vert}
\newcommand{\norm}[1]{\Vert #1 \Vert}

\newcommand{\fnorm}[1]{\norm{#1}_{\mathrm{F}}}

\newcommand{\minimize}{\text{minimize}}

\usepackage{xspace}\clearpage{}

\usepackage{lipsum}
\usepackage{amssymb,amsfonts}
\usepackage{graphicx}
\usepackage{epstopdf}
\ifpdf
  \DeclareGraphicsExtensions{.eps,.pdf,.png,.jpg}
\else
  \DeclareGraphicsExtensions{.eps}
\fi

\newsiamremark{remark}{Remark}
\newsiamremark{hypothesis}{Hypothesis}
\crefname{hypothesis}{Hypothesis}{Hypotheses}
\newsiamthm{claim}{Claim}
\newsiamthm{fact}{Fact}

\usepackage[noend]{myalgpseudocode}
\usepackage{algorithmicx}

\headers{Fast Randomized Subspace Algorithms}{Y.~Nakatsukasa and J.~A.~Tropp}

\title{Fast \& Accurate Randomized Algorithms \\ for Linear Systems and Eigenvalue Problems\thanks{Date: 29 October 2021.  Revised: 4 February 2022.
\funding{
JAT acknowledges ONR BRC N00014-1-18-2363 and NSF FRG 1952777.}}}

\author{Yuji Nakatsukasa\thanks{Mathematical Institute, University of Oxford, Oxford, UK 
  (\email{nakatsukasa@maths.oxford.ac.uk}).}
\and Joel A.~Tropp\thanks{Computing and Mathematical Sciences, Caltech, Pasadena, CA, USA 
  (\email{jtropp@cms.caltech.edu}).}}

\usepackage{amsopn}

\ifpdf
\hypersetup{
  pdftitle={Fast Randomized Algorithms for Linear Systems and Eigenvalue Problems},
  pdfauthor={Y. Nakatsukasa and J. A. Tropp}
}
\fi

\newcommand{\jatnote}[1]{\textcolor{blue}{\textbf{[JAT: #1]}}}
\newcommand{\ynnote}[1]{\textcolor{green}{\textbf{[YN: #1]}}}
\usepackage[normalem]{ulem}
\newcommand{\ignore}[1]{}
\DeclareMathOperator*{\subjectto}{subject\ to}

\begin{document}

\maketitle

\begin{abstract}
This paper develops a new class of algorithms for general linear systems and eigenvalue problems.
These algorithms apply fast randomized sketching to accelerate
subspace projection methods, such as GMRES and Rayleigh--Ritz.
This approach offers great flexibility in designing the basis
for the approximation subspace, which can improve scalability in many
computational environments.
As compared with classic methods, the new algorithms have
similar accuracy but run much faster.
For model problems, numerical experiments
show large advantages over MATLAB's optimized routines, including a
$100 \times$ speedup over \texttt{gmres} and a $10 \times$ speedup over
\texttt{eigs}. 
\end{abstract}



\section{Introduction}

Arguably, the most exciting recent development in numerical linear algebra (NLA)
is the advent of new randomized algorithms that are fast, scalable, robust, and reliable. For example, many practitioners have adopted  the ``randomized SVD'' and its relatives~\cite{halko2011finding,MartinssonTroppacta}
to compute truncated singular value decompositions of large matrices.
Randomized preconditioning~\cite{rokhlin2008fast,avron2010blendenpik}
allows us to solve highly overdetermined least-squares problems faster
than any previous algorithm.

In spite of these successes, our community has made less progress on other core challenges
from NLA, especially problems involving nonsymmetric square matrices.  
This paper exposes a new class of algorithms for solving general linear systems
and eigenvalue problems.  Our framework combines subspace projection
methods \cite{Saa03:Iterative-Methods,Saa11:Numerical-Methods}, such as GMRES and the Rayleigh--Ritz process,
with the modern technique of randomized sketching~\cite{sarlos2006improved,woodruff2014sketching,MartinssonTroppacta}.
This approach allows us to accelerate the existing methods by incorporating approximation subspaces
that are easier to construct.
The resulting algorithms are faster than their classic counterparts,
without much loss of accuracy.   In retrospect, the marriage of these ideas
appears inevitable.

\subsection{Sketching a least-squares problem}
\label{sec:sketch-intro}

The \textit{sketch-and-solve} paradigm~\cite{sarlos2006improved,woodruff2014sketching,MartinssonTroppacta}
is a basic tool for randomized matrix computations.
The idea is to decrease the dimension of a large problem by projecting it
onto a random subspace and to solve the smaller problem instead.
The solution of this ``sketched problem'' sometimes serves
in place of the solution to the original computational problem.

For a typical example, consider the $n \times d$ overdetermined least-squares problem
\begin{equation} \label{eqn:ls-intro}
\minimize_{\vct{y} \in \CC^d} \quad \norm{ \mtx{M}\vct{y} - \vct{f} }_2,
\end{equation}
where $\mtx{M} \in \CC^{n \times d}$ is a tall matrix with $n \gg d$.
The right-hand side $\vct{f} \in \CC^n$ and $\norm{\cdot}_2$ denotes the $\ell_2$ norm.
Draw a random sketching matrix $\mtx{S} \in \CC^{s \times n}$
with embedding dimension $s = 2d$, say.  Then solve the smaller $s \times d$ sketched problem
\begin{equation} \label{eqn:ls-sketch-intro}
\minimize_{\vct{y} \in \CC^d} \quad \norm{ \mtx{S} (\mtx{M} \vct{y} - \vct{f} ) }_2.
\end{equation}
For a carefully designed, ``fast'' sketching matrix $\mtx{S}$, the whole sketch-and-solve process may be
significantly faster than solving~\cref{eqn:ls-intro} directly. See~\cref{sec:subspace-embedding} for details.

We can compare the residual norms of the solution $\hat{\vct{y}}$
to the sketched problem~\cref{eqn:ls-sketch-intro} and the
solution $\vct{y}_{\star}$ to the original problem~\cref{eqn:ls-intro}.
The sketching method ensures that
\begin{equation} \label{eqn:ls-residuals-intro}
\norm{ \mtx{M} \vct{y}_{\star} - \vct{f}}_2 \quad\leq\quad
\norm{ \mtx{M} \hat{\vct{y}} - \vct{f} }_2 \quad\leq\quad 6 \cdot \norm{ \mtx{M} \vct{y}_{\star} - \vct{f}}_2.
\end{equation}
\textit{Provided that the original problem has a tiny residual}, the solution to the sketched
problem also yields a tiny residual!

\subsection{Solving linear systems by sketched GMRES}
\label{sec:linsys-intro}

Now, suppose that we wish to solve the (nonsymmetric, nonsingular) linear system
\begin{equation} \label{eqn:linsys-intro}
\text{Find $\vct{x} \in \CC^n$ : } \quad
\mtx{A} \vct{x} = \vct{f}
\quad\text{where $\mtx{A} \in \CC^{n \times n}$ and $\vct{f} \in \CC^n$.}
\end{equation}
All algorithms in this paper access the matrix via
products: $\vct{x} \mapsto \mtx{A} \vct{x}$.
Our approach builds on a standard template, called a \textit{subspace projection method}
\cite{Saa03:Iterative-Methods},
which casts the linear system as a variational problem.  We can treat
this formulation by sketching.  Let us summarize the ideas;
a full exposition appears in~\cref{sec:gmres,sec:gmres-basis}.

\subsubsection{Sketched GMRES}
\label{sec:sGMRES-intro}

For the moment, suppose that we have acquired a tall matrix $\mtx{B} \in \CC^{n \times d}$,
called a \textit{basis}, with the property that $\range(\mtx{B})$ contains a good approximate
solution to the linear system~\cref{eqn:linsys-intro}. That is, $\mtx{AB}\vct{y} \approx \vct{f}$
for some $\vct{y} \in \CC^d$. In addition, assume we have
the reduced matrix $\mtx{AB} \in \CC^{n \times d}$ at hand.
In typical situations, the basis has very low dimension: $d \ll n$. 

At its heart, the GMRES algorithm~\cite{gmres,SS85:Conjugate-Gradient-Like,Saa03:Iterative-Methods}
is a subspace projection method that
replaces the linear system~\cref{eqn:linsys-intro}
with the overdetermined least-squares problem
\begin{equation} \label{eqn:linsys-proj-intro}
\minimize_{\vct{y} \in \CC^d}\quad \norm{ \mtx{AB} \vct{y} - \vct{f} }_2.
\end{equation}
The solution $\vct{y}_{\star}$ to~\cref{eqn:linsys-proj-intro} yields an approximate solution $\vct{x}_{\mtx{B}} = \mtx{B}\vct{y}_{\star}$
to the linear system~\cref{eqn:linsys-intro}.  The residual norm $\norm{ \mtx{A} \vct{x}_{\mtx{B}} - \vct{f} }_2$
reflects how well the basis $\mtx{B}$ captures a solution to the linear system.

The least-squares formulation~\cref{eqn:linsys-proj-intro} is a natural candidate for sketching.
Draw a sketching matrix $\mtx{S} \in \CC^{s \times n}$ with $s = 2d$, say, and sketch the problem:
\begin{equation} \label{eqn:linsys-proj-sketch-intro}
\minimize_{\vct{y} \in \CC^d}\quad \norm{ \mtx{S}(\mtx{AB} \vct{y} - \vct{f}) }_2.
\end{equation}
The solution $\hat{\vct{y}}$ of the sketched problem~\cref{eqn:linsys-proj-sketch-intro}
induces an approximate solution $\hat{\vct{x}} = \mtx{B}\hat{\vct{y}}$ to the
linear system~\cref{eqn:linsys-intro}.  According to~\cref{eqn:ls-residuals-intro},
the residual norm $\norm{ \mtx{A} \hat{\vct{x}} - \vct{f} }_2$
is within a constant factor of the original residual norm $\norm{ \mtx{A} \vct{x}_{\mtx{B}} - \vct{f} }_2$. In summary, the sketched formulation~\cref{eqn:linsys-proj-sketch-intro}
is effective if and only if the subspace $\range(\mtx{B})$
contains an accurate approximate solution of the linear system.

We refer to~\cref{eqn:linsys-proj-sketch-intro} as the \textit{sketched GMRES} problem (sGMRES).
For an unstructured basis $\mtx{B}$, the sGMRES approach is faster than solving the original least-squares problem~\cref{eqn:linsys-proj-intro},
both in theory and in practice.
With careful implementation, sGMRES is reliable and robust, even when
the conditioning of the reduced matrix $\mtx{AB}$ is poor.  Indeed, it suffices
that $\kappa_2(\mtx{AB}) \lesssim u^{-1}$ where $u$ is the unit roundoff.\footnote{In standard IEEE double-precision arithmetic, the unit roundoff $u\approx 10^{-16}$.}
As a consequence, we have an enormous amount of flexibility in choosing the basis $\mtx{B}$.

\subsubsection{Krylov subspaces}

To make sGMRES work well, we must construct a subspace that captures an approximate solution to the linear system~\cref{eqn:linsys-intro}.
To that end, consider a Krylov subspace of the form
\begin{equation} \label{eqn:krylov-intro}
\set{K}_p(\mtx{A}, \vct{f})
	:= \lspan\{ \vct{f}, \mtx{A} \vct{f}, \mtx{A}^2 \vct{f}, \dots, \mtx{A}^{p-1} \vct{f} \}.
\end{equation}
The Krylov subspace often contains an excellent
approximate solution to the linear system, even
when the depth $p \ll n$.  See~\cite[Chaps.~6 and 7]{Saa03:Iterative-Methods}.

For computations, we need an explicit basis $\mtx{B}$ whose columns span
the Krylov subspace. Although it is straightforward to form
the monomial basis visible in~\cref{eqn:krylov-intro},
the condition number may grow exponentially~\cite{beckermann2000condition,Gau79:Condition-Polynomials},
rendering the basis useless for numerical purposes.
Instead, we will consider other procedures that quickly
construct Krylov subspace bases with smaller condition number.
\Cref{sec:gmres-basis} outlines several possible approaches.

For concreteness, we focus on the \textit{$k$-truncated Arnoldi process};
the parameter $k$ is a small natural number.
This algorithm assembles a basis $\mtx{B} = [\vct{b}_1, \dots, \vct{b}_d] \in \CC^{n \times d}$
iteratively.  Define $\vct{b}_{-i} = \vct{0}$ for $i \geq 0$.  Set $\vct{b}_1 = \vct{f} / \norm{\vct{f}}_2$.
For each $j = 2, \dots, d$, \begin{equation}
  \label{eq:partorth}
\vct{b}_j = \vct{w}_j / \norm{\vct{w}_j}_2
\quad\text{where}\quad
\vct{w}_j = (\Id - \vct{b}_{j-1} \vct{b}_{j-1}^* - \dots - \vct{b}_{j-k} \vct{b}_{j-k}^*) (\mtx{A}\vct{b}_{j-1}).  
\end{equation}
We have written ${}^*$ for the (conjugate) transpose.
Note that we obtain the reduced matrix $\mtx{AB}$ as a by-product of this computation.

We have found that $k$-truncated Arnoldi is often effective, even with $k = 2$ or $k = 4$.
Nevertheless, we are not aware of any fast, universal procedure
for constructing a Krylov subspace basis with full numerical rank,
short of strategies that perform costly full orthogonalization.  This is
a matter for further research.

\begin{algorithm}[t]\footnotesize
  \caption{sGMRES + $k$-truncated Arnoldi
  \label{alg:sgmres}}
  \begin{algorithmic}[1]
  \vspace{0.5pc}

	\Require{Matrix $\mtx{A} \in \CC^{n\times n}$, right-hand side $\vct{f} \in \CC^n$, initial guess $\vct{x} \in \CC^n$, basis dimension $d$, number $k$ of vectors for truncated orthogonalization, stability tolerance $\texttt{tol}=O(u^{-1})$.}

\Ensure{Approximate solution $\hat{\vct{x}} \in \CC^n$ to linear system~\cref{eqn:linsys-intro} and estimated
    residual norm $\hat{r}_{\mathrm{est}}$}

\vspace{0.5pc}

	\Function{sGMRES}{}

	\State	Draw subspace embedding $\mtx{S} \in \CC^{s \times n}$ with $s = 2(d + 1)$
		\Comment{See \cref{sec:random-embedding}} 

	\State  Form residual and sketch: $\vct{r} = \vct{f} - \mtx{A}\vct{x}$ and $\vct{g} = \mtx{S}\vct{r}$
		
	\State	Normalize basis vector $\vct{b}_1 = \vct{r} / \norm{\vct{r}}_2$ and apply matrix $\vct{m}_1 = \mtx{A}\vct{b}_1$

	\For{$j = 2, 3, 4, \dots, d$}
	
		\State	Truncated Arnoldi:
		$\mtx{w}_j = (\Id - \vct{b}_{j-1} \vct{b}_{j-1}^* - \dots - \vct{b}_{j-k} \vct{b}_{j-k}^*)\vct{m}_{j-1}$
		\Comment{$\vct{b}_{-i} = \vct{0}$ for $i \geq 0$}
		
		\State	 Normalize basis vector $\vct{b}_j = \vct{w}_j / \norm{\vct{w}_j}_2$ and apply matrix $\vct{m}_j = \mtx{A}\vct{b}_j$		
	\EndFor
		
		\State Sketch reduced matrix: $\mtx{C} = \mtx{S} [\vct{m}_1, \dots, \vct{m}_d]$
		
		\State Thin \textsf{QR} factorization: $\mtx{C} = \mtx{UT}$
		
		\If{condition number $\kappa_2(\mtx{T}) > \texttt{tol}$} warning...
		\State	Either whiten $\mtx{B} \gets \mtx{B} \mtx{T}^{-1}$ or form new residual and restart
		\Comment{See \cref{sec:adaptive-restart}}
		\EndIf
		
		\State	Solve least-squares problem: $\hat{\vct{y}} = \mtx{T}^{-1} (\mtx{U}^*\vct{g})$
		\Comment{See \cref{eqn:sgmres-yhat}}
	
		\State 	Residual estimate: $\hat{r}_{\mathrm{est}} = \norm{ (\Id - \mtx{UU}^*) \vct{g} }_2$
		\Comment{See \cref{eqn:sgmres-rhat}}
		
		\State	Construct solution: $\hat{\vct{x}} = \vct{x} + [\vct{m}_1,\dots,\vct{m}_j]\hat{\vct{y}}$
	\EndFunction

  \vspace{0.25pc}
	
	\Implementation{In line 6, use double Gram--Schmidt for stability.  In line 9, the \textsf{QR} factorization may require pivoting. In lines 11--12, apply $\mtx{T}^{-1}$ via triangular substitution.}

  \vspace{0.25pc}
\end{algorithmic}
\end{algorithm}

\begin{figure}[t] \includegraphics[width=0.505\textwidth]{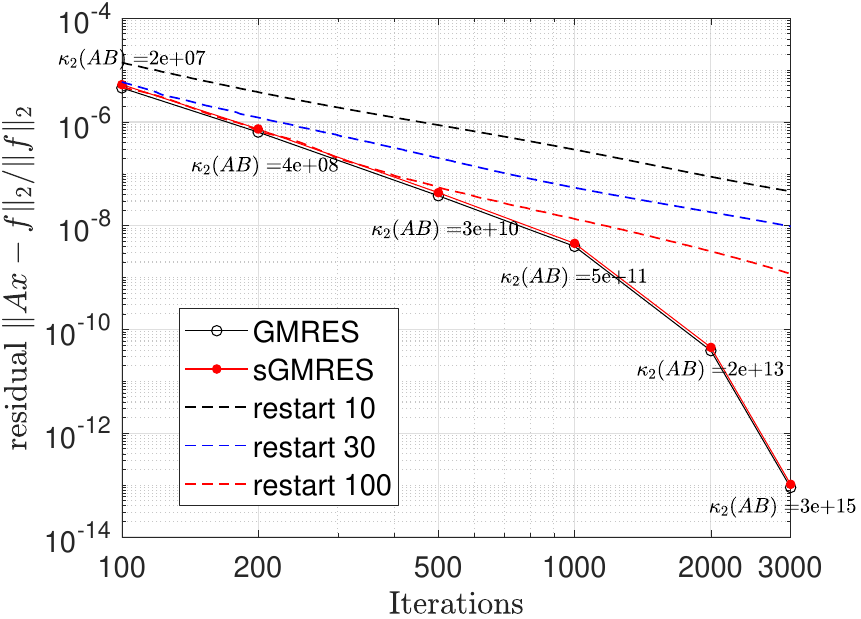}  
\includegraphics[width=0.48\textwidth]{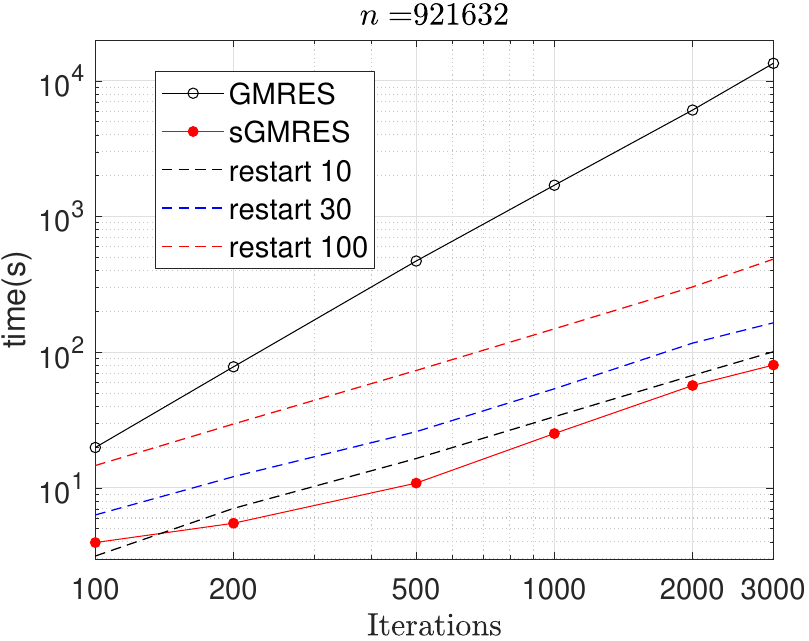}  
\caption{\textbf{GMRES versus sGMRES: Nonsymmetric linear system.}  These panels compare the performance of
MATLAB \texttt{gmres} (with and without restarting) against the sGMRES algorithm
(where the basis $\mtx{B}$ is computed by $k$-truncated Arnoldi with $k=4$).
The sparse linear system $\mtx{A}\vct{x} = \vct{f}$ has dimension $n=921,632$. 
\textbf{Left:} Relative residual
and the condition number $\kappa_2(\mtx{AB})$ of the reduced matrix. \textbf{Right:} Total runtime including basis generation. 
}\label{fig:gmres-intro}
\end{figure}

\subsubsection{Comparison with GMRES}

The standard version of GMRES~\cite{gmres} applies the expensive Arnoldi
process (with full orthogonalization; see~\cref{sec:arnoldi}) to build an orthonormal basis for the Krylov subspace,
and it exploits the structure of this basis to solve
the least-squares problem~\cref{eqn:linsys-proj-intro} efficiently.

In contrast, we propose to use a quick-and-dirty construction,
such as the $k$-truncated Arnoldi process, to obtain a basis for the Krylov subspace.
Then we solve the sGMRES least-squares problem~\cref{eqn:linsys-proj-sketch-intro}
to produce an approximate solution of the linear system.
When the basis dimension $d \ll n$, the sGMRES approach has lower arithmetic costs than classic GMRES,
while attaining similar accuracy:
\begin{center} \vspace{0.5pc}
\framebox{
GMRES: $O(nd^2)$ operations
\quad vs. \quad
sGMRES: $O(d^3 + nd \log d)$ operations.}
\vspace{0.5pc}
 \end{center}
This expression assumes that sGMRES uses $k$-truncated Arnoldi for $k$ constant,
as well as a fast sketching matrix (\cref{sec:random-embedding}).
See~\cref{alg:sgmres} for pseudocode.
 
As evidence for the benefits of using sGMRES, \Cref{fig:gmres-intro} depicts an over \textbf{100$\times$ speedup}
for a sparse nonsymmetric linear system with dimension $n = 921,632$.
In this case, sGMRES with $4$-truncated Arnoldi attains the same accuracy as GMRES
with full orthogonalization. sGMRES is comparable in speed to restarted GMRES
with restarting frequency $10$, whose convergence is significantly impaired.
\Cref{sec:experiments} provides more details on the experimental setup, as well as further illustrations.
For example, when applied to a sparse positive-definite linear system,
sGMRES can produce $\ell_2$ residual norms about $5\times$ smaller than the conjugate gradient (CG) method
after the same running time. Thus, it can be argued that sGMRES combines the speed of CG with the generality and robustness of GMRES.

\subsection{Solving eigenvalue problems by sketched Rayleigh--Ritz}

Similar ideas apply to spectral computations.
We pose the nonsymmetric eigenvalue problem
\begin{equation} \label{eqn:eig-intro}
\quad\text{Find nonzero $\vct{x} \in \CC^n$ and $\lambda \in \CC$ : }\quad
\mtx{A} \vct{x} = \lambda \vct{x}
\quad\text{where $\mtx{A} \in \CC^{n \times n}$.}
\end{equation}
As before, we access the matrix via products: $\vct{x} \mapsto \mtx{A}\vct{x}$.
Typically, we seek a family of eigenvectors associated
with a particular class of eigenvalues (e.g., largest real part, closest
to zero).  Let us outline a sketched subspace projection method for the eigenvalue problem.
Full details appear in~\cref{sec:eigs,sec:eigs-basis}. 

\subsubsection{Sketched Rayleigh--Ritz}
\label{sec:sRR-intro}

As in \cref{sec:sGMRES-intro}, suppose that we have procured a basis
$\mtx{B} \in \CC^{n \times d}$ and the reduced matrix $\mtx{AB} \in \CC^{n \times d}$.
The range of the basis should contain approximate eigenpairs $(\vct{x}, \lambda)$ for which
$\mtx{A}\vct{x} \approx \lambda \vct{x}$.  In this setting,
the most commonly employed strategy is the Rayleigh--Ritz (RR) method.

We begin with the classic variational formulation~\cite[Thm.~11.4.2]{parlettsym}
of RR:
\begin{equation} \label{eqn:RR-intro}
\minimize_{\mtx{M} \in \CC^{d \times d}}\quad
\fnorm{ \mtx{AB} - \mtx{BM} }.
\end{equation}
The solution is $\mtx{M}_{\star} = \mtx{B}^\dagger \mtx{A} \mtx{B}$,
where the dagger ${}^\dagger$ denotes the Moore--Penrose pseudoinverse.
At this point, RR frames the $d \times d$ eigenvalue problem
$\mtx{M}_{\star} \vct{y} = \theta \vct{y}$.  Each solution yields an approximate
eigenpair $(\mtx{B}\vct{y}, \theta)$ of the matrix $\mtx{A}$.

Evidently, the least-squares problem~\cref{eqn:RR-intro} is ripe for sketching.
Draw a sketching matrix $\mtx{S} \in \CC^{s \times n}$ with $s = 4d$, say,
and pass to the sketched RR problem:
\begin{equation} \label{eqn:pre-sRR-intro}
\minimize_{\mtx{M} \in \CC^{d\times d}} \quad
\fnorm{ \mtx{S}(\mtx{AB} - \mtx{BM}) }.
\end{equation}
We can compute the solution $\hat{\mtx{M}} = (\mtx{SB})^\dagger(\mtx{SAB})$ to the sketched
problem~\cref{eqn:pre-sRR-intro} faster than we can
obtain $\mtx{M}_{\star}$.
As before, we frame an ordinary eigenvalue problem:
\begin{equation} \label{eqn:pre-sRR-eig}
\hat{\mtx{M}} \vct{y} = \theta \vct{y}.
\end{equation}
For each solution
$(\hat{\vct{y}}, \hat{\theta})$, we obtain an approximate
eigenpair $(\mtx{B} \hat{\vct{y}}, \hat{\theta})$ of the
original matrix $\mtx{A}$.  We will show---both theoretically
and empirically---that the computed eigenpairs of~\cref{eqn:pre-sRR-eig}
are competitive with the eigenpairs obtained from RR.

We refer to~\cref{eqn:pre-sRR-intro} as the \textit{sketched Rayleigh--Ritz} (sRR) formulation.
Although it demands a careful implementation, sRR is faster than the
original least-squares method~\cref{eqn:RR-intro} for an unstructured basis $\mtx{B}$.
Moreover, sRR is robust, even when the basis $\mtx{B}$ has poor conditioning.
Indeed, it suffices that $\kappa_2(\mtx{B}) \lesssim u^{-1}$.

\begin{algorithm}[t]\footnotesize
  \caption{sRR + $k$-truncated Arnoldi} 
  \label{alg:srr}
  \begin{algorithmic}[1]
  \vspace{0.5pc}

	\Require{Matrix $\mtx{A}\in\CC^{n \times n}$, initial vector $\vct{b} \in \CC^n$, basis dimension $d$, number $k$ of vector for partial orthogonalization, stability tolerance $\texttt{tol}=O(u^{-1})$, convergence tolerance $\tau$.}

\Ensure{Approximate eigenpairs 
$(\vct{x}_i, \lambda_i)$ such that $\mtx{A}\vct{x}_i \approx \lambda_i \vct{x}_i$ and estimated residual norms $\hat{r}_{\mathrm{est},i}$.}

\vspace{0.5pc}

	\Function{sRR}{}

	\State	Draw subspace embedding $\mtx{S} \in \CC^{s \times n}$ with $s = 4d$ \Comment{See \cref{sec:random-embedding}} 

	\State	Starting vector: $\vct{w}_1 = \texttt{randn}(n,1)$
	
	\State	Normalize basis vector $\vct{b}_1 = \vct{w}_1 / \norm{\vct{w}_1}_2$ and apply matrix $\vct{m}_1 = \mtx{A}\vct{b}_1$

	\For{$j = 2, 3, 4, \dots, d$}
	
		\State	Truncated Arnoldi:
		$\mtx{w}_j = (\Id - \vct{b}_{j-1} \vct{b}_{j-1}^* - \dots - \vct{b}_{j-k} \vct{b}_{j-k}^*)\vct{m}_{j-1}$
		\Comment{$\vct{b}_{-i} = \vct{0}$ for $i \geq 0$}
		
		\State	Normalize $\vct{b}_j = \vct{w}_j / \norm{\vct{w}_j}_2$ and apply matrix  $\vct{m}_j = \mtx{A}\vct{b}_j$

	\EndFor
		
		\State Sketch basis
$\mtx{C} = \mtx{S}[\vct{b}_1,\ldots,\vct{b}_{d_{\max}}]$ and reduced matrix
$\mtx{D} = \mtx{S}[\vct{m}_1,\ldots,\vct{m}_{d_{\max}}]$
		\State	Thin \textsf{QR} factorization: $\mtx{C} = \mtx{UT}$
		
		\If{$\kappa_2(\mtx{T}) > \texttt{tol}$} warning:
			\State	Either whiten $\mtx{B}\leftarrow \mtx{B}\mtx{T}^{-1}$ or stabilize and solve~\cref{eq:qz} \Comment{See~\cref{sec:srr-stabilization}}
\EndIf
		
		\State	Solve eigenvalue problem: $\mtx{T}^{-1}\mtx{U}^*\mtx{D}\vct{y}_i=\lambda_i\vct{y}_i$ for $i=1,\ldots, d$ 
		\Comment{See \cref{eqn:Mhat-form}}
			
\State  Form residual estimates $\|\mtx{D}\vct{y}_i-\lambda_i\mtx{C}\vct{y}_i\|_2/\|\mtx{C}\vct{y}_i\|_2$
				\Comment{See \cref{eqn:srr-resid-est}, \cref{sec:srr-implementation}}

		\State	Identify set $\mathcal{I}$ of indices $i$ where residual is at most $\tau$
		\State Compute $\vct{x}_i=\mtx{B}\vct{y}_i$ and normalize $\vct{x}_i:=\vct{x}_i/\|\vct{x}_i\|_2$ for $i\in\mathcal{I}$, and output $(\vct{x}_i, \lambda_i)$

	\EndFunction

  \vspace{0.25pc}

	\Implementation{In line 6, use double Gram--Schmidt for stability.  In line 9, the \textsf{QR} factorization may require pivoting. In lines 11--12, apply $\mtx{T}^{-1}$ via triangular substitution.}
  \vspace{0.25pc}

\end{algorithmic}
\end{algorithm}

\subsubsection{Comparison with Arnoldi + Rayleigh--Ritz}

As before, we can deploy the Krylov subspace~\cref{eqn:krylov-intro}
for eigenvalue computations~\cite{parlettsym,Saa11:Numerical-Methods}.
In this case, we typically use a \textit{random} starting vector $\vct{\omega} \in \CC^n$
to generate the subspace $\set{K}_p(\mtx{A}; \vct{\omega})$.

To solve a large nonsymmetric eigenvalue problem,
one standard algorithm \cite[Sec.~6.2]{Saa11:Numerical-Methods}
applies the Arnoldi process (with full orthogonalization; see~\cref{sec:arnoldi})
to form an orthonormal basis for the Krylov subspace, and it uses the structure of
the basis to solve the RR eigenvalue problem efficiently.

Instead, we propose to combine a fast construction of a
Krylov subspace basis, such as $k$-truncated Arnoldi~\cref{eq:partorth},
with the sRR eigenvalue problem~\cref{eqn:pre-sRR-eig}.
When the basis dimension $d \ll n$, this algorithm uses less arithmetic
than the classic approach:
\begin{center} \vspace{0.5pc}
\framebox{
RR: $O(nd^2)$ operations
\quad vs. \quad
sRR: $O(d^3 + nd \log d)$ operations.}
\vspace{0.5pc}
 \end{center}
This expression includes basis generation via $k$-truncated Arnoldi for $k$ constant,
and sRR uses a fast sketching matrix (\cref{sec:random-embedding}).
See~\cref{alg:srr} for pseudocode.

\begin{figure}[t] \includegraphics[width=0.50\textwidth]{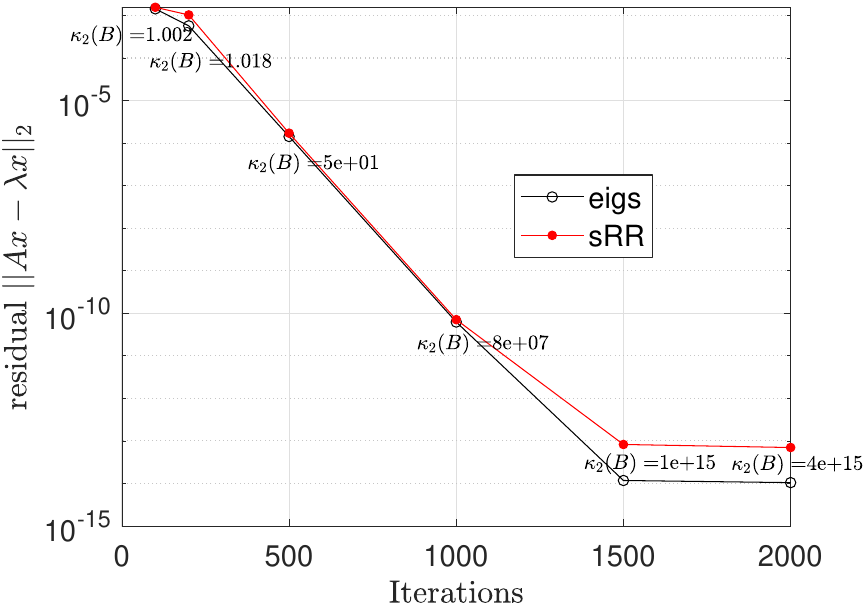}  
\includegraphics[width=0.47\textwidth]{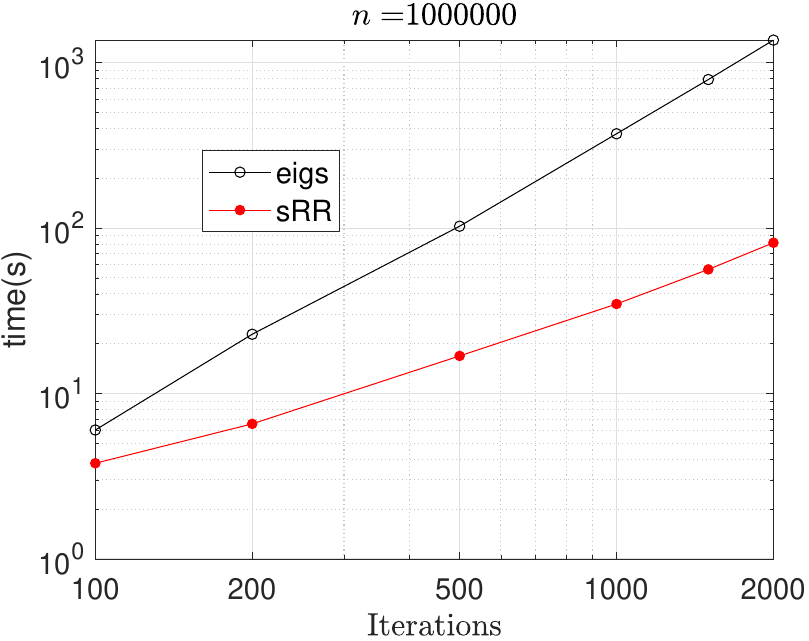}  
\caption{\textbf{RR versus sRR: Nonsymmetric eigenvalue problem.}  These panels compare the performance of MATLAB \texttt{eigs}
against the sRR algorithm (where the basis $\mtx{B}$ is computed by $k$-truncated Arnoldi with $k = 10$).
The sparse eigenvalue problem $\mtx{A}\vct{x} = \lambda \vct{x}$
has dimension $n = 10^6$, and it arises from a trust-region subproblem in optimization.
\textbf{Left:} Relative residual for the right-most eigenpair and the condition number $\kappa_2(\mtx{B})$ of the basis. 
Throughout, all eigenvectors are normalized to have unit norm.
\textbf{Right:} Total runtime including basis generation.
}\label{fig:eigs-intro}
\end{figure}

As evidence, \Cref{fig:eigs-intro} highlights
an eigenvalue problem arising in numerical optimization
for which sRR runs over \textbf{10$\times$ faster} than the MATLAB
\texttt{eigs} command.  Even so, both methods
compute the desired eigenpair to the same accuracy.
\Cref{sec:experiments} describes the experimental setup and provides
further illustrations.  For example, when applied to a sparse symmetric
eigenvalue problem, sRR can outperform standard implementations of
the Lanczos method in both speed and reliability.

\subsubsection{Block Krylov subspaces}

For eigenvalue problems, there is also a compelling opportunity to explore
alternative subspace constructions.
For example, consider the \textit{block} Krylov subspace
\begin{equation} \label{eqn:rand-Krylov-intro}
\set{K}_p(\mtx{A}, \mtx{\Omega})
	:=\lspan\{ \mtx{\Omega}, \mtx{A\Omega}, \mtx{A}^2 \mtx{\Omega}, \dots, \mtx{A}^{p-1} \mtx{\Omega} \}
\quad\text{where $\mtx{\Omega} \in \CC^{n \times b}$.}
\end{equation}
We commonly generate the Krylov subspace from a random
matrix $\mtx{\Omega}$.  The standard prescription
recommends a small block size $b$ and a large depth $p$,
but recent research~\cite[Sec.~11]{MartinssonTroppacta} has
shown the value of a large block size $b$ and a small depth $p$.

We must take care in constructing the block Krylov subspace.
Truncated Arnoldi is only competitive when the block size $b$ is
a small constant. For larger $b$, the Chebyshev recurrence offers an elegant
way to form a basis $\mtx{B} = [\mtx{B}_1, \dots, \mtx{B}_p] \in \CC^{n \times bp}$:
$$
\mtx{B}_1 = \mtx{\Omega}; \qquad
\mtx{B}_2 = \mtx{A\Omega}; \qquad
\mtx{B}_i = 2 \mtx{A}\mtx{B}_{i-1} - \mtx{B}_{i-2}
\quad\text{for $i = 3, \dots, p$.}
$$
We obtain the reduced matrix $\mtx{AB}$ as a by-product.
In practice, the Chebyshev polynomials must be shifted and scaled to adapt
to the spectrum of $\mtx{A}$.  See \cref{sec:eigs-basis} for details and
alternative methods for fast basis construction.

\subsection{Discussion}

Our idea to combine subspace projection methods with sketching
offers compelling advantages over the classic algorithms,
especially in modern computing environments.
Nevertheless, it must be acknowledged that this approach suffers
from some of the same weaknesses as GMRES and RR.  For example,
when the basis $\mtx{B}$ is a Krylov subspace, these methods
are limited by the approximation power of Krylov subspaces.
Furthermore, we are not aware of a universal method for quickly
computing a full-rank basis for the Krylov subspace, short of
expensive strategies based on full orthogonalization.
Both of these points merit further attention.

With hindsight, our framework appears as an obvious application
of the sketch-and-solve paradigm for overdetermined least-squares
problems.  A critical reader may even wonder whether this idea is
actually novel.  Let us respond to this concern.

We believe that we are the first authors to identify the natural
connection between sketching and subspace projection methods for
linear algebra problems.  Indeed, we are not aware of any
prior work where authors sketch the minimum-residual
formulations~\cref{eqn:linsys-proj-sketch-intro,eqn:pre-sRR-intro} of
general linear systems and eigenvalue problems.

Second, to obtain algorithms that are asymptotically faster than
classic methods, we must also employ efficient constructions of Krylov
subspace bases.  In the past, researchers have regarded these
techniques as a way to \textit{postpone} expensive orthogonalization
steps in parallel computing environments~\cite{JC91:Parallelizable-Restarted,philippe2012generation}.
In contrast, sketching sometimes allows us to \textit{eliminate}
the orthogonalization steps.  Thus, we can finally take full
advantage of the potential of fast computational bases.

In particular, there is a remarkable opportunity to design a new
class of preconditioners for iterative solution of large-scale
linear algebra problems.  Indeed, since our approach allows for
more flexible bases and mitigates orthogonalization costs, we can
still derive benefits from a mediocre preconditioner that only
reduces the iteration complexity to 100s or 1000s of iterations.
We are excited about this prospect.

\subsection{Related work}

Within the NLA literature, there has been relatively little research on the application of sketching for high-accuracy
matrix computations.  After this paper was completed, we learned about several papers that involve both subspace
projection methods and sketching, but they are different in spirit from our work.

First, Balabanov \& Nouy~\cite{BN19:Randomized-Linear,BN21:Randomized-Linear-II} have argued that randomized algorithms
may be used to accelerate reduced-order modeling.  They focus on sketching an approximation subspace that captures solutions
to a parameterized system of linear equations.   In this setting, many additional complications arise from the need to scope out
the parameter manifold and to interpolate between solutions at different parameter values.  Their work does not suggest
the simple, fast algorithms that we propose here.

Second, Balabanov \& Grigori~\cite{balabanov2020randomized} have observed that sketching can reduce the
cost of orthogonalization in the (block) Arnoldi method.  Although this is an elegant idea, it only reduces the
leading constant in the cost of the basis construction, but not its asymptotic scaling.  In contrast, our
approach yields asymptotically faster algorithms.  As suggested in their follow-up work~\cite{BG21:Randomized-Block},
there are appealing opportunities to combine their approach with ours.

There is also extensive work on sketching in the theoretical algorithms literature and machine learning literature;
for example, see the surveys~\cite{woodruff2014sketching,DM18:Lectures-Randomized,MartinssonTroppacta}.  For the most part,
the algorithms that follow from the traditional sketching perspective lead to low-accuracy output that would not be considered
acceptable by numerical analysts.

\subsection{Roadmap}

In~\cref{sec:subspace-embedding}, we give a rigorous treatment
of sketching for least-squares problems.  \Cref{sec:gmres,sec:gmres-basis,sec:sgmres-algs}
develop and analyze the sGMRES method and associated basis
constructions.  \Cref{sec:eigs,sec:eigs-basis} contain the
analogous developments for sRR.  Computational experiments
in~\cref{sec:experiments} confirm that these algorithms are
fast, robust, and reliable.  
\Cref{sec:extension,sec:prospects} describe extensions and prospects. 

\subsection{Notation}

The symbol ${}^*$ denotes the (conjugate) transpose of a vector or matrix.
We write $\norm{\cdot}_2$ for the $\ell_2$ norm or the spectral norm, while $\fnorm{\cdot}$ is the Frobenius norm.
The dagger $\dagger$ denotes the pseudoinverse.  For a matrix $\mtx{M} \in \CC^{n \times d}$,
define the largest singular value $\sigma_{\max}(\mtx{M}) := \sigma_1(\mtx{M})$ and the minimum singular value
$\sigma_{\min}(\mtx{M}) := \sigma_{\min\{n,d\}}(\mtx{M})$.  The condition number
$\kappa_2(\mtx{M}) := \sigma_{\max}(\mtx{M}) / \sigma_{\min}(\mtx{M})$.

\section{Background: Subspace embeddings}
\label{sec:subspace-embedding}

A subspace embedding is a linear map, usually from a high-dimensional space to a low-dimensional space,
that preserves the $\ell_2$ norm of every vector in a given subspace.
This definition is due to Sarl{\'o}s~\cite{sarlos2006improved};
see also~\cite{woodruff2014sketching,MartinssonTroppacta}. 

\begin{definition}[Subspace embedding]
Suppose that the columns of $\mtx{B} \in \CC^{n \times d}$ span the subspace $\set{L} \subseteq \CC^n$.
A matrix $\mtx{S} \in \CC^{s \times n}$ is called a subspace embedding for $\set{L}$ with
distortion $\eps \in (0, 1)$ if
\begin{equation} \label{eqn:subspace-embed}
(1 - \eps) \cdot \norm{ \mtx{B} \vct{y} }_2 \leq \norm{ \mtx{SB} \vct{y} }_2
	\leq (1+\eps) \cdot \norm{ \mtx{B} \vct{y}  }_2
	\quad\text{for all $\vct{y} \in \CC^d$.}
\end{equation}
\end{definition}

For matrix computations, we need to design subspace embeddings that have several additional properties.
First, the subspace embedding $\mtx{S}$ should be equipped with a fast matrix--vector
multiply so that we can perform the data reduction process efficiently.
Second, the subspace $\set{L}$ is typically unknown, so we must draw a subspace embedding
\textit{at random} to achieve~\cref{eqn:subspace-embed} with high probability.
Last, to randomly embed a $d$-dimensional subspace with distortion $\eps$,
the optimal scaling of the embedding dimension $s$ follows the law
$s \approx d / \eps^2$.  Owing to this relation, subspace embeddings are only appropriate
in settings where a moderate distortion, say $\eps = 1/\sqrt{2}$, is enough for computational purposes.

Before turning to constructions in~\cref{sec:random-embedding},
let us outline the applications of subspace embeddings
that we will use in this paper.

\subsection{Sketching for least-squares problems}

As discussed in~\cref{sec:sketch-intro},
we can use a subspace embedding to reduce the dimension of an overdetermined least-squares problem.
This idea is also due to Sarl{\'o}s~\cite{sarlos2006improved};
it serves as the foundation for a collection of
methods called the \textit{sketch-and-solve} paradigm~\cite{woodruff2014sketching,MartinssonTroppacta}.  

\begin{fact}[Sketching for least-squares] \label{fact:sketching}
Let $\mtx{M} \in \CC^{n \times d}$ be a matrix, and suppose that $\mtx{S} \in \CC^{s \times n}$
is a subspace embedding for $\range([\mtx{M}, \vct{f}])$ with distortion $\eps \in (0,1)$.
For every vector $\vct{y} \in \CC^d$, we have the two-sided inequality
\begin{equation} \label{eqn:ls-sketch}
(1 - \eps) \cdot \norm{  \mtx{M} \vct{y} - \vct{f} }_2
	\leq \norm{ \mtx{S}(\mtx{M} \vct{y} - \vct{f}) }_2
	\leq (1 + \eps) \cdot \norm{  \mtx{M} \vct{y} - \vct{f} }_2.
\end{equation}
In particular, the solution $\vct{y}_{\star}$ to the least-squares problem~\cref{eqn:ls-intro} and
the solution $\hat{\vct{y}}$ to the sketched least-squares problem~\cref{eqn:ls-sketch-intro} satisfy
residual norm bounds
\begin{equation} \label{eqn:ls-residuals}
\norm{ \mtx{M} \vct{y}_{\star} - \vct{f} }_2 \leq \norm{ \mtx{M} \hat{\vct{y}} - \vct{f} }_2
	\leq \frac{1 + \eps}{1 - \eps} \cdot \norm{ \mtx{M} \vct{y}_{\star} - \vct{f} }_2.
\end{equation}
Equation~\cref{eqn:ls-residuals} justifies the claim~\cref{eqn:ls-residuals-intro} with $\eps = 1/\sqrt{2}$.
\end{fact}

\subsection{Whitening the basis}

Rokhlin \& Tygert~\cite{rokhlin2008fast} observed that a subspace embedding yields an
inexpensive way to precondition an iterative algorithm for the overdetermined
least-squares problem.  We can invoke the same idea to approximately orthogonalize,
or \textit{whiten}, a given basis.

\begin{fact}[Whitening]
Let $\mtx{B} \in \CC^{n \times d}$ be a basis with full column rank.
Let $\mtx{S} \in \CC^{s \times n}$ be a subspace embedding for $\range(\mtx{B})$
with distortion $\eps \in (0, 1)$.  Compute a \textsf{QR} factorization of the
sketched basis: $\mtx{SB} = \mtx{UT}$ with $\mtx{U} \in \CC^{s \times d}$ orthonormal and $\mtx{T} \in \CC^{d\times d}$ permuted triangular.
Then the whitened basis $\bar{\mtx{B}} = \mtx{BT}^{-1}$ satisfies
\begin{equation} \label{eqn:whitening}
\kappa_2(\bar{\mtx{B}}) = \frac{\sigma_{\max}(\bar{\mtx{B}})}{\sigma_{\min }(\bar{\mtx{B}})}
	\leq \frac{1+\eps}{1 - \eps}.
\end{equation}
Furthermore, we have the condition number diagnostic
\begin{equation} \label{eqn:cond-diagnostic}
\frac{1-\eps}{1+\eps} \cdot \kappa_2(\mtx{T}) \leq \kappa_2(\mtx{B}) \leq \frac{1+\eps}{1-\eps} \cdot \kappa_2(\mtx{T}).
\end{equation}
\end{fact}

\subsection{Constructing a subspace embedding}
\label{sec:random-embedding}

There are many performant constructions of fast randomized subspace embeddings
that work for an unknown subspace of bounded dimension~\cite[Sec.~9]{MartinssonTroppacta}.
Let us summarize two that are most relevant for our purposes.  In each case,
for a subspace with dimension $d$,
to obtain empirical distortion $\eps \in (0 , 1)$,
we set the embedding dimension $s = d / \eps^2$.
We focus on the complex field; modifications
for the real field are straightforward.

\subsubsection{SRFTs}

First, we introduce the subsampled random Fourier transform (SRFT)~\cite{AC09:Fast-Johnson-Lindenstrauss,woolfe2008fast,tropp2011improved,MartinssonTroppacta}.
This subspace embedding\footnote{For worst-case problems, a more elaborate SRFT construction may be needed~\cite[Sec.~9]{MartinssonTroppacta}.}
 takes the form
\begin{equation} \label{eqn:SRFT}
\mtx{S} = \sqrt{\frac{n}{s}} \mtx{DFE} \in \CC^{s \times n}.
\end{equation}
In this expression, $\mtx{D} \in \CC^{s \times n}$ is a diagonal projector onto $s$ coordinates, chosen independently at random,
$\mtx{F} \in \CC^{n \times n}$ is the unitary discrete Fourier transform (DFT), and $\mtx{E} \in \CC^{n \times n}$
is a diagonal matrix whose entries are independent Steinhaus\footnote{A \textit{Steinhaus random variable} is uniform on the complex unit circle $\{ z \in \CC : \abs{z} = 1 \}$.}
random variables.  The cost of applying the matrix $\mtx{S}$ to an
$n\times d$ matrix is $O(nd \log d)$ operations using the subsampled FFT algorithm~\cite{woolfe2008fast}.

\subsubsection{Sparse maps}

Next, we describe the sparse dimension reduction map~\cite{MM13:Low-Distortion-Subspace,nelson2013osnap,clarkson2017low,cohen2016nearly,MartinssonTroppacta},
which is useful for sparse data and may require less data movement.
It takes the form
\begin{equation} \label{eqn:sparse-map}
\mtx{S} = \frac{1}{\sqrt{s}} [ \vct{s}_1, \dots, \vct{s}_n ] \in \CC^{s \times n}.
\end{equation}
The columns of $\mtx{S}$ are statistically independent.  Each column $\vct{s}_i$ has exactly $\zeta$ nonzero
entries, drawn from the Steinhaus distribution, placed in uniformly random coordinates.
For reliability, we choose the sparsity level $\zeta = \lceil 2 \log(1 + d) \rceil$.
We can apply $\mtx{S}$ to a matrix $\mtx{M}$ with $O(\zeta \cdot \texttt{nnz}(\mtx{M}))$ operations,
but it may require a sparse arithmetic library to achieve the best performance.

\section{Solving linear systems with sGMRES}
\label{sec:gmres}

We return to the linear system
\begin{equation} \label{eqn:linsys}
\text{Find $\vct{x} \in \CC^n$ : }\quad
\mtx{A}\vct{x} = \vct{f}
\quad\text{where $\mtx{A} \in \CC^{n \times n}$ and $\vct{f} \in \CC^n$.}
\end{equation}
This section elaborates on the sGMRES method outlined
in~\cref{sec:linsys-intro}.
\Cref{sec:gmres-basis} discusses methods
for constructing the basis required by sGMRES.
\Cref{sec:sgmres-algs} combines these ideas
to obtain complete sGMRES algorithms.

\subsection{Derivation of GMRES}

Fix a \textit{full rank} basis $\mtx{B} \in \CC^{n \times d}$
and the reduced matrix $\mtx{AB} \in \CC^{n \times d}$.
Suppose that $\vct{x}_0 \in \CC^n$ is an initial guess
for the solution of~\cref{eqn:linsys} with residual
$\vct{r}_0 := \vct{f} - \mtx{A}\vct{x}_0$.  Lacking
prior information, we may take $\vct{x}_0 = \vct{0}$.

Consider the affine family of approximate solutions
to~\cref{eqn:linsys}
of the form $\vct{x} = \vct{x}_0 + \mtx{B}\vct{y}$
where $\vct{y} \in \CC^d$.  Among this class, we may select a representative whose
residual $\vct{r} = \vct{f} - \mtx{A}\vct{x} = \vct{r}_0 - \mtx{AB}\vct{y}$
has the minimum $\ell_2$ norm:
\begin{equation} \label{eqn:gmres}
\minimize_{\vct{y} \in \CC^d}\quad
\norm{ \mtx{AB}\vct{y} - \vct{r}_0 }_2. 
\end{equation}
With some imprecision, we refer to~\cref{eqn:gmres} as the GMRES problem~\cite{gmres}. By calculus, the least-squares problem~\cref{eqn:gmres}
is equivalent to the normal equations:
\begin{equation} \label{eqn:gmres-orthogonality}
\text{Find $\vct{y} \in \CC^d$ :}\qquad
	(\mtx{AB})^*(\mtx{AB}\vct{y} - \vct{r}_0) = \mtx{0}.
\end{equation}
We can stably solve \cref{eqn:gmres} using a \textsf{QR} factorization of the reduced matrix $\mtx{AB}$~\cite[Ch.~20]{Higham:2002:ASNA}.
The cost is $O(nd^2)$ arithmetic operations, assuming that the reduced matrix $\mtx{AB}$ is unstructured.
Given a solution $\vct{y}_{\mtx{B}}$ to either problem, we obtain a new approximate solution $\vct{x}_{\mtx{B}} = \vct{x}_0 + \mtx{B}\vct{y}_{\mtx{B}}$
with residual $\vct{r}_{\mtx{B}} = \vct{f} - \mtx{A}\vct{x}_{\mtx{B}}$.

The formulation~\cref{eqn:gmres-orthogonality}
is also called a \textit{Petrov--Galerkin method}~\cite[Chap.~5]{Saa03:Iterative-Methods}
with approximation space $\vct{x}_0 + \range(\mtx{B})$
and orthogonality space $\range(\mtx{AB})$.
The GMRES algorithm~\cite{SS85:Conjugate-Gradient-Like,gmres} is a particular instance
where $\mtx{B}$ is an orthonormal basis
for a Krylov subspace generated by $\vct{r}_0$.
GMRES forms the basis $\mtx{B}$ via the Arnoldi process (\cref{sec:arnoldi}),
which involves $d$ matvecs with $\mtx{A}$ plus $O(nd^2)$ arithmetic.
This reduces~\cref{eqn:gmres} to a structured least-squares problem
that can be solved in $O(d^2)$ operations.

\subsection{Derivation and analysis of sGMRES}

To develop the sGMRES method, we just sketch the GMRES
problem~\cref{eqn:gmres}.
Construct a subspace embedding $\mtx{S} \in \CC^{s \times n}$
for $\range( [\mtx{AB}, \vct{r}_0] )$ with distortion $\eps \in (0,1)$.
The sketched GMRES problem is
\begin{equation} \label{eqn:sgmres}
\minimize_{\vct{y} \in \CC^d} \quad
\norm{ \mtx{S}(\mtx{AB} \vct{y} - \vct{r}_0) }_2.
\end{equation}
Let $\hat{\vct{y}} \in \CC^d$ denote the solution of~\cref{eqn:sgmres}.
Write $\hat{\vct{x}} = \vct{x}_0 + \mtx{B} \hat{\vct{y}}$
and $\hat{\vct{r}} = \vct{f} - \mtx{A}\hat{\vct{x}}$.

We have an \textit{a priori} comparison of the GMRES~\cref{eqn:gmres} and sGMRES~\cref{eqn:sgmres}
residual norms because of the relation~\cref{eqn:ls-residuals}:
\begin{equation} \label{eqn:sgmres-residual}
\norm{ \mtx{A}\vct{x}_{\mtx{B}} - \vct{f} }_2
	\leq \norm{ \mtx{A} \hat{\vct{x}} - \vct{f} }_2
	\leq \frac{1+\eps}{1-\eps} \cdot \norm{ \mtx{A} \vct{x}_{\mtx{B}} - \vct{f} }_2.
\end{equation}
Thus, sGMRES produces approximate solutions to~\cref{eqn:linsys} with small $\ell_2$ residuals
precisely when GMRES does.
\textit{A posteriori}, we can diagnose the quality of
the computed solution $\hat{\vct{x}}$ by examining
the sketched residual norm:
\begin{equation} \label{eqn:sgmres-residual-norm}
\hat{r}_{\mathrm{est}} := \norm{ \mtx{S} (\mtx{AB} \hat{\vct{y}} - \vct{r}_0) }_2
	\in [1 - \eps, 1 + \eps] \cdot \norm{ \mtx{A} \hat{\vct{x}} - \vct{f} }_2.
\end{equation}
The last display is a consequence of~\cref{eqn:ls-sketch}.

For both GMRES~\cref{eqn:gmres} and sGMRES~\cref{eqn:sgmres},
the fundamental challenge is to produce a basis $\mtx{B}$ that captures
an approximate solution to the linear system~\cref{eqn:linsys}.
We return to this matter in~\cref{sec:gmres-basis}.

\subsection{Implementation}
\label{sec:sgmres-implementation}

Let us outline a numerically robust implementation of sGMRES
and describe some of the issues that arise.

The algorithm operates with either an SRFT~\cref{eqn:SRFT} or a sparse embedding~\cref{eqn:sparse-map},
depending on which is more appropriate to the computational environment.
We recommend the embedding dimension $s = 2(d+1)$, which typically yields
distortion $\eps = 1/\sqrt{2}$.  In view of~\cref{eqn:sgmres-residual},
the sGMRES residual norm is less than $6 \times$
the GMRES residual norm, although the discrepancy is often smaller in practice.

To obtain the data for the sGMRES problem~\cref{eqn:sgmres},
we sketch the reduced matrix ($\mtx{SAB} \in \CC^{s \times d}$) and
the right-hand side ($\mtx{S}\vct{r}_0 \in \CC^s$)
at a cost of $O(nd \log d)$ operations.
To solve~\cref{eqn:sgmres},
we compute a thin, pivoted \textsf{QR} decomposition of the sketched matrix:
$\mtx{SAB} = \mtx{UT}$ where $\mtx{U} \in \CC^{s \times d}$ is orthonormal
and $\mtx{T} \in \CC^{d\times d}$ is a triangular matrix with permuted columns.
A minimizer of the sGMRES problem is
\begin{equation} \label{eqn:sgmres-yhat}
\hat{\vct{y}} = (\mtx{SAB})^\dagger (\mtx{S}\vct{r}_0)
	= \mtx{T}^{-1} (\mtx{U}^* (\mtx{Sr}_0)).
\end{equation}
Of course, we apply the inverse by triangular substitution.
The sketched residual norm~\cref{eqn:sgmres-residual-norm} admits the simple expression
\begin{equation} \label{eqn:sgmres-rhat}
\hat{r}_{\mathrm{est}} = \norm{ (\Id - \mtx{UU}^*) (\mtx{S}\vct{r}_0) }_2.
\end{equation}
The two preceding displays require $O(d^3)$ arithmetic since $s = O(d)$.
Last, we explicitly form the approximate solution
$\hat{\vct{x}} = \vct{x}_0 + \mtx{B} \hat{\vct{y}}$ at a cost of $O(nd)$ operations.

In summary, given the basis $\mtx{B}$, the cost of forming and solving the sGMRES problem~\cref{eqn:sgmres}
is $O(d^3 + nd \log d)$ arithmetic.
In contrast, for an unstructured basis, the cost of solving the GMRES problem~\cref{eqn:gmres} is $O(nd^2)$ arithmetic.
\Cref{sec:sgmres-algs} provides an accounting
of the costs of forming the basis and solving the least-squares problem.

\subsection{Stability}
\label{sec:sgmres-stabilization}

The classical stability result~\cite[Thm.~20.3]{Higham:2002:ASNA} shows that standard numerical methods for the least-squares problem~\cref{eqn:sgmres} produce a solution with essentially optimal residual as long as $\kappa_2(\mtx{SAB})\lesssim u^{-1}$. According to~\cref{eqn:cond-diagnostic}, this condition is equivalent to $\kappa_2(\mtx{AB})\lesssim u^{-1}$.

Our computational work (\cref{sec:experiments}) confirms that
sGMRES is reliable unless the reduced matrix $\mtx{AB}$ is very badly conditioned.  
In our experience, it suffices that $\kappa_2(\mtx{AB}) \leq 10^{14}$ in double-precision arithmetic.
Therefore, we have wide latitude to design bases that we can construct quickly;
see~\cref{sec:gmres-basis}.
We will provide evidence that sGMRES with a fast basis construction
is more efficient than GMRES with a structured basis.

\subsection{Restarting}

Standard implementations of GMRES periodically restart~\cite[Sec.~6.5.5]{Saa03:Iterative-Methods}.
That is, they use a basis $\mtx{B}$ to compute an approximate solution $\vct{x}_{\mtx{B}}$
to the linear system~\cref{eqn:linsys} with the residual vector
$\vct{r}_{\mtx{B}} = \vct{r}_0 - \mtx{A} \vct{x}_{\mtx{B}}$.
If the residual norm $\norm{\vct{r}_{\mtx{B}}}_2$ exceeds an error tolerance,
the residual vector $\vct{r}_{\mtx{B}}$ is used to generate a new basis,
which is fed back to GMRES to construct another approximate
solution.  This process is repeated until a solution of desired quality is  obtained.

Restarting has a number of benefits for the process of basis construction.
It allows us to work with bases that have fewer columns, which limits the
cost of storing the basis.  For orthogonal basis constructions,
restarting reduces the cost of orthogonalization.  For non-orthogonal
basis constructions, the restarting process helps control the conditioning
of the basis.  On the other hand, restarted GMRES may not converge if
the bases are not rich enough (see \cref{fig:gmres-intro}), and we pay
for the convergence delay with every restart.
As we will discuss in~\cref{sec:sgmres-algs},
sGMRES can help us manage all of these issues.

\subsection{Preconditioning}

For difficult linear systems, we may need a preconditioner
$\mtx{P} \in \CC^{n \times n}$ to solve it successfully
with either GMRES or sGMRES.
The preconditioned system has the form
\begin{equation} \label{eqn:precond}
\mtx{P}^{-1} \mtx{A} \vct{x} = \mtx{P}^{-1} \vct{f}.
\end{equation}
A good preconditioner has two features~\cite[Chaps.~9 and 10]{Saa03:Iterative-Methods}.  First,
the matrix $\mtx{P}^{-1}\mtx{A}$
has a more ``favorable'' structure than $\mtx{A}$.
Second, we can solve $\mtx{P}\vct{z} = \vct{g}$
efficiently.  (Let us emphasize that we only interact with
$\mtx{P}^{-1}$ by solving linear systems!)
Although preconditioning is
critical in practice, it is heavily problem dependent,
so we will not delve into examples.

We may derive sGMRES for the preconditioned system~\cref{eqn:precond},
following the same pattern as before.  Note that we employ the preconditioned matrix
$\mtx{P}^{-1} \mtx{A}$ when we construct the basis $\mtx{B}$ and
the reduced matrix $\mtx{P}^{-1} (\mtx{AB})$.
The details are routine.  We believe that
sGMRES opens up new opportunities for
designing preconditioners because it
is faster and more flexible than GMRES.

\section{Constructing a basis for sGMRES}
\label{sec:gmres-basis}

As we have seen, the success of both GMRES~\cref{eqn:gmres} and sGMRES~\cref{eqn:sgmres}
hinges on the approximation power of the basis.
Krylov subspaces are, perhaps, the most natural way
to capture solutions to a linear system when we access
the matrix via products~\cite[Chaps.~6 and 7]{Saa03:Iterative-Methods}.
In this section, we describe
a number of ways to compute non-orthogonal bases for Krylov subspaces.
Although these strategies are decades old, they warrant
a fresh look because sGMRES has a fundamentally different
computational profile from GMRES.

\subsection{The single-vector Krylov subspace}

Many iterative methods for solving the linear system~\cref{eqn:linsys}
implicitly search for solutions in the Krylov subspace
$$
\set{K}_{p}(\mtx{A}; \vct{r})
	:= \lspan\{ \vct{r}, \mtx{A}\vct{r}, \mtx{A}^2 \vct{r}, \dots, \mtx{A}^{p-1} \vct{r} \}
	= \lspan\{ \varphi(\mtx{A}) \vct{r} : \deg(\varphi) \leq p-1 \}.
$$
In this context, the generating vector $\vct{r} \in \CC^{n}$ is often the normalized
residual $\vct{r}_0 / \norm{\vct{r}_0}_2$, defined by $\vct{r}_0 = \vct{f} - \mtx{A}\vct{x}_0$,
where $\vct{x}_0$ is an approximate solution to~\cref{eqn:linsys}.
The function $\varphi$ ranges over polynomials with degree at most $p-1$.

A basis $\mtx{B} \in \CC^{n \times d}$ for the Krylov subspace
$\set{K}_p(\mtx{A}; \vct{r})$ comprises a system of vectors
that spans the subspace.  We can write
$$
\mtx{B} = [\vct{b}_1, \dots, \vct{b}_d]
\quad\text{where}\quad
\vct{b}_j = \varphi_j(\mtx{A}) \vct{r}
\quad\text{for $j =1, \dots,d$.}
$$
The filter polynomials $(\varphi_j : j = 1,\dots, d)$ have degree
at most $p - 1$, and they are usually linearly independent
(so $d = p$).  In most cases, the polynomials are also graded
$(\deg(\varphi_j) = j - 1)$, and they are constructed sequentially by a recurrence.
This process delivers the reduced matrix $\mtx{AB}$ without any extra work.

For example, the monomial basis takes the form $\vct{b}_1 = \vct{r}$ and $\vct{b}_j = \mtx{A} \vct{b}_{j-1}$
for $j = 2, \dots, p$.  The associated polynomials are
$\varphi_j(t) = t^{j-1}$ for $j = 1, \dots, p$.  For many matrices $\mtx{A}$,
the conditioning of the monomial basis for $\set{K}_p(\mtx{A}; \vct{r})$
grows exponentially with $p$, so it is inimical to numerical computation~\cite{Gau79:Condition-Polynomials}.

We will consider other constructions of Krylov subspace bases that are more suitable in practice.  Our aim is to control the resources
used to obtain the basis, including arithmetic, (working) storage, communication,
synchronization, etc.   We can advance these goals by relaxing
the requirement that the basis be well-conditioned.

For theoretical analysis of the approximation power of Krylov subspaces
in the context of linear system solvers, see~\cite[Sec.~6.11]{Saa03:Iterative-Methods}.

\subsection{The Arnoldi process}
\label{sec:arnoldi}

It is supremely natural to build an \textit{orthonormal} basis
$\mtx{Q} \in \CC^{n \times p}$ for the Krylov subspace
$\set{K}_p(\mtx{A}; \vct{r})$ sequentially.  This is called
the Arnoldi process~\cite[Sec.~6.3]{Saa03:Iterative-Methods}.
The initial vector $\vct{q}_1 = \vct{r} /\norm{\vct{r}}_2$.
After $j$ steps, the method updates the partial basis
$\mtx{Q}_j = [ \vct{q}_1, \dots, \vct{q}_j ]$ by appending the vector
$$
\vct{q}_{j+1} = \vct{w}_{j+1} / \norm{\vct{w}_{j+1}}_2
\quad\text{where}\quad
\vct{w}_{j+1} = (\Id - \mtx{Q}_j \mtx{Q}_j^*) (\mtx{A}\vct{q}_{j}).
$$
The Arnoldi basis $\mtx{Q}_p \in \CC^{n \times p}$ has the happy property that
$$
\mtx{AQ}_p = \mtx{Q}_p \mtx{H}_p + \vct{w}_{p} \mathbf{e}_{p}^*
\quad\text{where $\mtx{H}_p \in \CC^{p \times p}$ is upper Hessenberg.}
$$
As a consequence, we can solve the least-squares problem~\cref{eqn:gmres}
with $\mtx{B} = \mtx{Q}_p$ in $O(p^2)$ time and produce the approximate
solution $\vct{x}_{\mtx{B}}$ in $O(np)$ operations.  This is roughly how the
standard implementation of the GMRES algorithm operates~\cite{gmres}.

The orthogonalization steps in the Arnoldi process are expensive.
For $p$ iterations, they expend $O(n p^2)$ arithmetic,
and they may also involve burdensome inner-products, communication,
and synchronization.
Robust implementations usually incorporate modified or double Gram--Schmidt
or else use Householder reflectors. 

The literature contains many strategies for controlling the
orthogonalization costs in the Arnoldi process~\cite[Chap.~6]{Saa03:Iterative-Methods}.
sGMRES motivates us to reevaluate techniques for building a nonorthogonal basis.
For example, we can use $k$-truncated Arnoldi, as in~\cref{eq:partorth},
which reduces the cost of basis generation to $O(n k^2)$.
Provided the reduced matrix $\mtx{AB}$ is reasonably conditioned,
we can still obtain accurate solutions to the linear system via
sGMRES~\cref{eqn:sgmres}.

Relatedly, Balabanov \& Grigori~\cite{balabanov2020randomized} have proposed
to use a low-dimensional sketch of the basis vectors to implement
an approximate orthogonalization strategy.  Their approach has the same
asymptotic cost as full orthogonalization, but similar ideas might be
invoked to accelerate partial or selective orthogonalization.

\subsection{The Lanczos recurrence}

For this subsection, assume $\mtx{A}$ is Hermitian.
In this case, the Arnoldi process simplifies to a three-term recurrence~\cite[Sec.~6.6]{Saa03:Iterative-Methods}:
$$
\vct{q}_{j+1} = \vct{w}_{j+1} / \norm{\vct{w}_{j+1}}_2
\quad\text{where}\quad
\vct{w}_{j+1} = (\Id - \vct{q}_{j}\vct{q}_j^* - \vct{q}_{j-1}\vct{q}_{j-1}^*) (\mtx{A}\vct{q}_j).
$$
The Lanczos basis $\mtx{Q}_p = [\vct{q}_1, \dots, \vct{q}_p] \in \CC^{n \times p}$
has the remarkable property that
\begin{equation}  \label{eq:lancos}
\mtx{AQ}_p = \mtx{Q}_p \mtx{J}_p + \vct{w}_p \vct{e}_p^*
\quad\text{where $\mtx{J}_p \in \CC^{p \times p}$ is tridiagonal.}  
\end{equation}
This allows us to solve the least-squares problem~\cref{eqn:gmres} with
$\mtx{B} = \mtx{Q}_p$ in $O(p)$ time, and we construct the approximate
solution $\vct{x}_{\mtx{B}}$ with $O(np)$ arithmetic.
This is roughly how the MINRES algorithm operates~\cite{minresoriginal}.

For $p$ iterations, the Lanczos recurrence costs just $O(np)$ operations,
but it has complicated behavior in finite-precision arithmetic.
This issue is not devastating when Lanczos is used to solve linear systems~\cite[Chap.~5]{LS12:Krylov-Subspace},
but it can present a more serious challenge when solving eigenvalue problems~\cite[Chap.~13]{parlettsym}.

Although it is very efficient to solve the least-squares problem~\cref{eqn:gmres}
by passing to the tridiagonal matrix $\mtx{J}_p$,
it is more reliable to sketch $\mtx{S}(\mtx{AQ}_p)$ and to solve the sketched
problem~\cref{eqn:sgmres} instead.  The approach based on sketching is
competitive with MINRES when $p \ll n$.

The literature describes many approaches for maintaining the orthogonality
of the Lanczos basis, such as selective orthogonalization~\cite[Chap.~13]{parlettsym}.
When using the basis for sGMRES, we might simply omit the extra orthogonalization steps.
Alternatively, we may adopt the approximate orthogonalization strategies sketched
in \cref{sec:arnoldi}.

\subsection{The Chebyshev recurrence}
\label{sec:chebyshev}

In some settings, we may wish to avoid the orthogonalization steps entirely
because they involve operations on high-dimensional basis vectors.
We can achieve this goal by using other polynomial recurrences to construct
a Krylov subspace basis. This idea is attributed to Joubert \& Carey~\cite{JC91:Parallelizable-Restarted}.

For simplicity, suppose that the spectrum of $\mtx{A}$ is contained in
the axis-aligned rectangle $[c \pm \delta_x, \pm \delta_y]$, and set $\varrho = \max\{\delta_x, \delta_y\}$.  Then we can assemble a shifted-and-scaled
Chebyshev basis $\mtx{B} \in \CC^{n\times p}$ via the following recurrence~\cite{Man77:Tchebychev-Iteration,JC91:Parallelizable-Restarted}.
Let $\vct{b}_1 = \vct{r}/\norm{\vct{r}}_2$ and $\vct{b}_2 = (2\varrho)^{-1} (\mtx{A} - c \Id) \vct{b}_1$.
Then
\begin{equation}
  \label{eq:Chebrecurse}
\vct{b}_j = \frac{1}{\varrho} \left[(\mtx{A} - c\Id) \vct{b}_{j-1} - \frac{\delta_x^2 - \delta_y^2}{4\varrho} \vct{b}_{j-2} \right]
\quad\text{for $j = 3, \dots, p$.}  
\end{equation}
In practice, we also rescale each basis vector $\vct{b}_j$ to have unit $\ell_2$ norm
after it has played its role in the recurrence.
The key theoretical fact is that the Chebyshev basis tends to have
a condition number that grows \textit{polynomially} in $p$, rather than exponentially.
This claim depends on assumptions that the eigenvalues of the matrix are
equidistributed over an ellipse~\cite{Gau72:Condition-Orthogonal,Man77:Tchebychev-Iteration,JC91:Parallelizable-Restarted,philippe2012generation}.

To implement this procedure, we may first apply a few iterations of the Arnoldi
method (\cref{sec:arnoldi-method}) to estimate the spectrum of $\mtx{A}$.  More generally, we find a (transformed)
ellipse that contains the spectrum.  Then we adapt the Chebyshev polynomials
to this ellipse~\cite{philippe2012generation}.  The overall cost of
constructing a Chebyshev basis for $\set{K}_p(\mtx{A}; \vct{r})$ 
is $O(np)$, and it involves no orthogonalization whatsoever.

\subsection{Newton polynomials}

The Newton polynomials provide another standard construction of
a nonorthogonal basis for the Krylov subspace~\cite{philippe2012generation}.
Suppose that $\theta_1, \dots, \theta_p \in \CC$ are complex-valued
shift parameters.  Then we can build a basis $\mtx{B} \in \CC^{n \times p}$
for $\set{K}_p(\set{A}; \vct{r})$ via the recurrence
$$
\vct{b}_1 = \vct{r} / \norm{\vct{r}}_2
\quad\text{and}\quad
\vct{b}_j = (\mtx{A} - \theta_{j-1} \Id)\vct{b}_{j-1}
\quad\text{for $j = 2, \dots, p$.}
$$
The shifts $\theta_i$ are often chosen to be estimated eigenvalues
of $\mtx{A}$, obtained from an invocation of the Arnoldi method (\cref{sec:arnoldi-method}).
The overall computational profile of constructing the Newton basis
is similar to constructing a Chebyshev basis.

\subsection{Local orthogonalization}

We can improve the conditioning of a computed basis
$\mtx{B} \in \CC^{n \times d}$ by local orthogonalization.
Indeed, it is generally helpful to orthogonalize subcollections
of basis vectors, even if it proves too expensive to orthogonalize
all of the basis vectors.  In particular, scaling each column
to have unit $\ell_2$ norm is always appropriate.
See~\cite{vdS69:Condition-Numbers} for an analysis.

\section{sGMRES algorithms}
\label{sec:sgmres-algs}

This section presents complete algorithms for solving the
linear system~\cref{eqn:linsys} via sGMRES, including options for
adaptive basis generation.

\subsection{Basic implementation}

\cref{alg:sgmres} contains pseudocode for a basic
implementation of sGMRES using the $k$-truncated
Arnoldi basis~\cref{eq:partorth}.  We recommend this version of the
algorithm when the user lacks
information about the spectrum of $\mtx{A}$.
Given bounds on the spectrum, one may
replace the truncated Arnoldi basis with
a Chebyshev basis (\cref{sec:chebyshev}).
\Cref{tab:sgmres-cost} summarizes the arithmetic costs.

\begin{table}[t]
\caption{\textbf{GMRES versus sGMRES: Arithmetic.}  This table compares the total arithmetic cost of solving an $n \times n$
linear system using a $d$-dimensional basis via standard GMRES and via sGMRES.  For sGMRES, we consider
both $k$-truncated Arnoldi and the Chebyshev basis.
Heuristically, the parameters $k \ll d \ll n$.  Constant factors are suppressed.} \label{tab:sgmres-cost}
\begin{center} \renewcommand{\arraystretch}{1.2}
\begin{tabular}{|l||c|c|c|c|c|c}
\hline
			& Matrix access				& Form basis	& Sketch		& LS solve	& Form soln.	\\
\hline
Std.~GMRES	& $d T_{\mathrm{matvec}}$	& $nd^2$		& ---			& $d^2$		& $nd$			\\
sGMRES-$k$	& $d T_{\mathrm{matvec}}$	& $ndk$			& $nd \log d$	& $d^3$ 	& $nd$			\\
sGMRES-Cheb & $d T_{\mathrm{matvec}}$	& $nd$			& $nd \log d$	& $d^3$		& $nd$			\\
\hline
\end{tabular}
\end{center}
\end{table}

\subsection{Iterative sGMRES}
\label{sec:sgmres-iterative}

As noted, 
most methods for producing the Krylov subspace basis
are recursive.  They generate the columns of $\mtx{B}$
and the reduced matrix $\mtx{AB}$ in sequence.
This observation suggests an iterative implementation
of sGMRES.  We sketch the columns of the
reduced matrix as they arrive, incrementally solving the sGMRES
problem~\cref{eqn:sgmres} at each step.

Let $d_{\max}$ be a user-specified parameter that bounds
the maximum depth allowed for the Krylov subspace.
Draw and fix a randomized subspace embedding
$\mtx{S} \in \CC^{s\times n}$ with embedding dimension
$s = 2 (d_{\max} + 1)$.  As we compute each column $\mtx{A} \vct{b}_j$ of the reduced matrix,
we immediately form the sketch $\mtx{S}(\mtx{A}\vct{b}_j)$
and update the \textsf{QR} decomposition:
\begin{equation} \label{eqn:sgmres-iterative-qr}
\mtx{S}(\mtx{A}\mtx{B}_j) = \mtx{U}_j \mtx{T}_j
\quad\text{where}\quad
\mtx{B}_j = [\vct{b}_1, \dots, \vct{b}_j].
\end{equation}
At each step, we obtain an approximate solution to the linear system:
\begin{equation} \label{eqn:sgmres-iterative-soln}
\hat{\vct{y}}_j = \mtx{T}_j^{-1} (\mtx{U}_j^* (\mtx{S}\vct{r}_0))
\quad\text{with}\quad
\hat{r}_{\mathrm{est},j} = \norm{ (\Id - \mtx{U}_j \mtx{U}_j)^* (\mtx{S}\vct{r}_0) }_2.
\end{equation}
Repeat this process until the estimated residual norm $\hat{r}_{\mathrm{est},j}$
is sufficiently small or we breach the threshold $d_{\max}$
for the size of the Krylov space.  After $d$ iterations,
the arithmetic costs of~\cref{eqn:sgmres-iterative-qr,eqn:sgmres-iterative-soln}
match the non-sequential implementation (\cref{sec:sgmres-implementation}) with a basis of size $d$.

\subsection{Adaptive restarting}
\label{sec:adaptive-restart}

There is a further opportunity to design an adaptive strategy for restarting.
According to~\cref{eqn:cond-diagnostic}, the condition number
$\kappa_2(\mtx{T}_j)$ is comparable with $\kappa_2(\mtx{AB}_j)$.  When first $\kappa_2(\mtx{T}_j) > \texttt{tol}$,
we recognize that it is time to restart.  We generate the new Krylov subspace using
the \textit{previous} residual $\hat{\vct{r}}_{j-1} = \vct{r}_0 - \mtx{AB}_{j-1} \hat{\vct{y}}_{j-1}$.
Alternatively, instead of restarting, we could approximately orthogonalize $\mtx{B}$ by replacing it with $\mtx{B}\mtx{T}^{-1}$, whose condition number is constant,
and continue generating basis vectors.

\subsection{Storage-efficient versions}

In situations where storage is at a premium, we can even avoid storing
the reduced matrix $\mtx{AB}_j$ by sketching its columns sequentially
and discarding them immediately after sketching.  Once the estimated
residual norm $\hat{r}_{\mathrm{est},j}$ is sufficiently small, we
can construct the approximate solution
$$
\hat{\vct{x}}_j = \vct{x}_0 + \mtx{AB}_j\hat{\vct{y}}_j
	= \vct{x}_0 + \sum\nolimits_{i=1}^j (\mtx{AB}_j)_i (\hat{\vct{y}}_j)_i
$$
by iteratively regenerating the columns of the reduced matrix $\mtx{AB}_j$
and linearly combining them on the fly.  For some basis constructions (e.g., truncated Arnoldi or Chebyshev),
we only need to maintain a few columns of $\mtx{B}$ and the $j$ columns of $\mtx{S}(\mtx{AB}_j)$.
This modification doubles the arithmetic cost associated with basis generation
(matvecs plus orthogonalization). A similar technique was used in~\cite{YTFUC21:Scalable-Semidefinite}.

\subsection{Obtaining a solution with full accuracy}

While the constant-factor loss~\cref{eqn:sgmres-residual} in sGMRES is unlikely to be an issue,
we can obtain a solution with the same quality as GMRES by using
$\mtx{T}$ as a preconditioner to solve~\cref{eqn:gmres} via
an iterative method as in~\cite{rokhlin2008fast,avron2010blendenpik}.
This method still requires $\kappa_2(\mtx{AB}) \lesssim u^{-1}$ to operate
reliably.

\section{The sketched Rayleigh--Ritz method}
\label{sec:eigs}

Let us turn to the nonsymmetric eigenvalue problem
\begin{equation} \label{eqn:eig}
\text{Find nonzero $\vct{x} \in \CC^n$ and $\lambda \in \CC$ : }\quad
\mtx{A}\vct{x} = \lambda \vct{x}
\quad\text{where $\mtx{A} \in \CC^{n \times n}$.}
\end{equation}
We will provide an implementation and analysis of the
sRR method outlined in \cref{sec:sRR-intro}.  \Cref{sec:symeig} describes modifications
for the symmetric eigenvalue problem.  \Cref{sec:eigs-basis} covers techniques
for constructing the basis for sRR.

\subsection{Perspectives on Rayleigh--Ritz}
\label{sec:rr-faces}

Fix a \textit{full-rank} basis $\mtx{B} \in \CC^{n \times d}$,
and let $\mtx{AB} \in \CC^{n \times d}$ be the reduced matrix.
Rayleigh--Ritz is best understood as a \textit{Galerkin method}
for computing eigenvalues~\cite[Sec.~4.3]{Saa11:Numerical-Methods}.
Among nonzero vectors of the form $\vct{x} = \mtx{By}$, we seek a residual $\vct{r} = \mtx{A}\vct{x} - \theta\vct{x}$
orthogonal to $\range(\mtx{B})$.  More precisely,
\begin{equation} \label{eqn:RR-PG}
\text{Find nonzero $\vct{y} \in \CC^d$ and $\theta \in \CC$ :} \qquad
\mtx{B}^* (\mtx{AB} \vct{y} - \theta \mtx{B} \vct{y}) = \mtx{0}.
\end{equation}
Rearranging, we see that~\cref{eqn:RR-PG} can be posed as an
ordinary eigenvalue problem:
\begin{equation} \label{eqn:RR-eig}
\mtx{M}_{\star} \vct{y} := \mtx{B}^\dagger (\mtx{AB}) \vct{y} = \theta \vct{y}
\quad\text{where $\vct{y} \neq \vct{0}$ and $\theta \in \CC$.}
\end{equation}
Recall that eigenvalue problems are invariant under similarity transforms.
In the present context, the computed eigenpairs only depend on the range of $\mtx{B}$,
so they are invariant under the map $\mtx{B}\leftarrow \mtx{BT}$
for a nonsingular $\mtx{T} \in \CC^{d \times d}$.
Therefore, if $\mtx{Q} \in \CC^{n \times d}$ is an orthonormal basis for
$\range(\mtx{B})$, then we may pass to \begin{equation} \label{eqn:RR-eig-ortho}
\mtx{Q}^* (\mtx{A} \mtx{Q}) \vct{z} = \theta \vct{z}
\quad\text{where $\vct{z} \neq \vct{0}$.}
\end{equation}
Given a solution $(\vct{z}, \theta)$ to~\cref{eqn:RR-eig-ortho}, we obtain an approximate eigenpair
$(\mtx{Q} \vct{z}, \theta)$ of the matrix $\mtx{A}$. 
This is the most typical presentation of RR.

In contrast, consider the problem of minimizing the residual
over the subspace:
\begin{equation} \label{eqn:rect-eig}
\minimize_{\vct{y} \in \CC^d, \theta \in \CC}\quad
\norm{ \mtx{AB}\vct{y} - \theta \mtx{B}\vct{y} }_2
\quad\text{subject to $\norm{ \mtx{B} \vct{y} }_2 = 1$.}
\end{equation}
This formulation is sometimes called a \textit{rectangular eigenvalue problem}~\cite{boutry2005generalized,itomurota2016}.
Let us emphasize that the RR method~\cref{eqn:RR-eig} \textit{does not} solve the rectangular
eigenvalue problem.  Nevertheless, for any eigenpair $(\vct{y}_{\star}, \theta_{\star})$
of the matrix $\mtx{M}_{\star}$, it holds that
\begin{equation} \label{eqn:RR-resid}
\norm{ \mtx{AB} \vct{y}_{\star} - \theta_{\star} \mtx{B} \vct{y}_{\star} }_2
	= \norm{ (\mtx{AB} - \mtx{BM}_{\star})\vct{y}_{\star} }_2.
\end{equation}
The matrix $\mtx{M}_{\star}$ from~\cref{eqn:RR-eig}
\textit{does} solve a related variational problem~\cite[Thm.~11.4.2]{parlettsym}:
\begin{equation} \label{eqn:RR-matrix}
\minimize_{\mtx{M} \in \CC^{d\times d}} \quad
\fnorm{ \mtx{AB} - \mtx{BM} }.
\end{equation}
These connections support the design and analysis of a sketched version of RR.

\subsection{The Arnoldi method}
\label{sec:arnoldi-method}

The Arnoldi method is a classic algorithm~\cite[Sec.~6.2]{Saa11:Numerical-Methods} for eigenvalue problems based on RR.
First, it invokes the Arnoldi process (\cref{sec:arnoldi}) to build an orthonormal
basis $\mtx{Q} \in \CC^{n \times d}$ for a Krylov subspace (generated by a random vector)
at a cost of $O(nd^2)$ operations.  This construction ensures that
$\mtx{Q}^*\mtx{AQ}$ has upper Hessenberg form,
so we can solve the eigenvalue problem~\cref{eqn:RR-eig} with
$O(d^2)$ operations by means of the QR algorithm~\cite[Chap.~7]{Saa11:Numerical-Methods}.  Each eigenpair $(\vct{y}, \theta)$ of~\cref{eqn:RR-eig}
induces an approximate eigenpair $(\mtx{B}\vct{y}, \theta)$
of $\mtx{A}$, which we can form with $O(nd)$ operations.

\subsection{Derivation of sRR}

We can view the sRR method as a sketched version of the matrix optimization
problem~\cref{eqn:RR-matrix}.  Consider a subspace embedding
$\mtx{S} \in \CC^{s \times n}$ for $\range([ \mtx{AB}, \mtx{B} ])$
with distortion $\eps \in (0,1)$.  The sketched problem is
\begin{equation} \label{eqn:sRR-matrix-form}
\minimize_{\mtx{M} \in \CC^{d\times d}} \quad \fnorm{ \mtx{S}(\mtx{AB} - \mtx{BM}) }.
\end{equation}
The sRR method finds a solution $\hat{\mtx{M}} \in \CC^{d \times d}$
to this optimization problem.  Then it poses the ordinary eigenvalue problem
\begin{equation} \label{eqn:srr-eig}
\hat{\mtx{M}} \vct{y} = \theta \vct{y}
\quad\text{where $\vct{y} \neq \vct{0}$.}
\end{equation}
This computation yields up to $d$ eigenpairs $(\hat{\vct{y}}_i, \hat{\theta}_i)$
of the matrix $\hat{\mtx{M}}$.  We obtain approximate eigenpairs of $\mtx{A}$
by the transformation $(\mtx{B}\hat{\vct{y}}_i, \hat{\theta}_i)$.

Sketching allows us to obtain inexpensive \textit{a posteriori} error bounds.
For a computed eigenpair $(\hat{\vct{y}}, \hat{\theta})$ of $\hat{\mtx{M}}$,
it is cheap to form the sketched residual:
\begin{equation} \label{eqn:srr-resid-est}
\hat{r}_{\mathrm{est}}
	:= \hat{r}_{\mathrm{est}}(\hat{\vct{y}}, \hat{\theta})
	:= \norm{ \mtx{S}(\mtx{AB}\hat{\vct{y}} - \hat{\theta}\mtx{B}\hat{\vct{y}})}_2/\norm{\mtx{SB}\hat{\vct{y}}}_2.
\end{equation}
By definition, the subspace embedding $\mtx{S}$ ensures that the true residual satisfies
\begin{equation}  \label{eq:rescompare}
\frac{1-\eps}{ 1+\eps} \cdot \hat{r}_{\mathrm{est}} 
\leq 
\frac{\norm{ \mtx{AB}\hat{\vct{y}} - \hat{\theta} \mtx{B}\hat{\vct{y}} }_2}{\norm{\mtx{B}\hat{\vct{y}}}_2}
\leq 
\frac{ 1+\eps}{1-\eps} \cdot \hat{r}_{\mathrm{est}}.  
\end{equation}
In other words, we can diagnose when the sRR method has (or has not)
produced a high-quality approximate eigenpair $(\mtx{B}\hat{\vct{y}}, \hat{\theta})$
of the original matrix $\mtx{A}$.

\subsection{Implementation of sRR}
\label{sec:srr-implementation}

To implement sRR, we may use either an SRFT embedding~\cref{eqn:SRFT}
or a sparse embedding~\cref{eqn:sparse-map}.  We recommend the
embedding dimension $s = 4d$, which typically results
in distortion $\eps = 1/\sqrt{2}$ for the range of $[\mtx{AB}, \mtx{B}]$.

We first sketch the reduced matrix ($\mtx{SAB} \in \CC^{s\times d}$)
and the basis ($\mtx{SB} \in \CC^{s\times d}$) at a cost of
$O(nd \log d)$ operations.  Next, we compute a thin, pivoted \textsf{QR}
decomposition $\mtx{SB} = \mtx{UT}$ of the sketched basis.
A minimizer of the sRR problem~\cref{eqn:sRR-matrix-form} is the matrix
\begin{equation} \label{eqn:Mhat-form}
\hat{\mtx{M}} := (\mtx{SB})^\dagger (\mtx{SAB})
	= \mtx{T}^{-1} (\mtx{U}^* (\mtx{SAB})) \in \CC^{d \times d}.
\end{equation}
We apply the inverse by triangular substitution.
Then invoke the QR algorithm to solve the
eigenvalue problem~\cref{eqn:srr-eig}.
Each of the last three steps costs $O(d^3)$ operations.

Given a computed eigenpair $(\hat{\vct{y}}, \hat{\theta})$,
we can obtain the sketched residual value
$\hat{r}_{\mathrm{est}}(\hat{\vct{y}}, \hat{\theta})$
from~\cref{eqn:srr-resid-est} at a cost of $O(d^2)$
operations.  If the residual estimate is sufficiently small,
we declare that $(\mtx{B}\hat{\vct{y}}, \hat{\theta})$
is an approximate eigenpair of $\mtx{A}$.
For maximum efficiency, we present the approximate
eigenvector $\hat{\vct{x}} = \mtx{B}\hat{\vct{y}} \in \CC^n$
in factored form.  If we need the full vector
$\hat{\vct{x}}$, it costs $O(nd)$ operations.
Ironically, if we extract a large number of explicit eigenvectors,
this last step dominates the cost of the computation.
Usually, the number of high-quality approximate eigenpairs is moderate.

In summary, given the basis $\mtx{B}$,
if we use sRR to solve~\cref{eqn:eig}, the cost of reporting
the factored form of $d$ approximate eigenpairs
is $O(d^3 + nd \log d)$ operations.  In contrast,
RR requires $O(nd^2)$ arithmetic with an unstructured
basis.
Our numerical experience indicates that sRR is a robust
alternative to RR so long as the condition number
of the basis $\kappa_2(\mtx{B}) \leq 10^{14}$.
This fact allows us to exploit fast non-orthogonal basis
constructions; see~\cref{sec:eigs-basis}.
See~\cref{alg:srr} for a simple implementation
of sRR with a partial Arnoldi basis.

\subsection{Stabilization}
\label{sec:srr-stabilization}

The output of sRR is almost identical to RR provided that $\kappa_2(\mtx{B}) \lesssim u^{-1}$.  This condition is very generous.
In contrast, recall that the standard stability analysis~\cite[Chap.~13]{parlettsym}
for the Lanczos algorithm asks that $\kappa_2(\mtx{B}) < 1 + \sqrt{u}$.

If we see that the sketched basis $\mtx{SB}$
is very badly conditioned ($\kappa_2(\mtx{SB})\gtrsim u^{-1}$),
then the condition number diagnostic~\cref{eqn:cond-diagnostic} implies
that the basis $\mtx{B}$ is also very badly conditioned.
In this case, we can stabilize sRR by regularizing the basis.
For example, we may compute the truncated SVD:
$$
\mtx{SB} = \mtx{U\Sigma V}^* + \mtx{E}
\quad\text{where 
$\norm{\mtx{\Sigma}}_2/\norm{\mtx{E}}_2 \leq \texttt{tol}$.}
$$
Then we use the QZ algorithm~\cite[\S 7.7]{golubbook4th} to solve the \textit{generalized} eigenvalue problem\footnote{When $\kappa_2(\mtx{B}) \gtrsim u^{-1}$, numerical experiments suggest this approach is more stable than reducing to a standard eigenvalue problem as in~\cref{eqn:Mhat-form}.}
\begin{equation}  \label{eq:qz}
(\mtx{U}^*\mtx{SAB})\mtx{V} \vct{z} = \theta \mtx{\Sigma} \vct{z}.  
\end{equation}
Each solution yields an sRR eigenpair $(\mtx{V}\vct{z}, \theta)$
and an associated approximate eigenpair $(\mtx{BV}\vct{z}, \theta)$
of $\mtx{A}$.  
The asymptotic cost is the same as the basic implementation.

\subsection{Why does sRR work?\nopunct}\label{sec:whysRR}

We will argue that RR and sRR solve eigenvalue problems
that are similar to a pair of nearby eigenvalue problems.

First, recall that $\mtx{SB} = \mtx{UT}$ is the \textsf{QR} decomposition
of the sketched basis.  Per~\cref{eqn:whitening}, the whitened
basis $\bar{\mtx{B}} := \mtx{BT}^{-1}$ has conditioning $\kappa_2(\bar{\mtx{B}}) \leq (1+\eps)/(1-\eps)$.
Now, consider the variational problem~\cref{eqn:RR-matrix}
with respect to the whitened basis $\bar{\mtx{B}}$
and its sketched version:
$$
\begin{aligned}
&\minimize_{\mtx{M}}\quad \fnorm{ \mtx{A}\bar{\mtx{B}} - \bar{\mtx{B}} \mtx{M} }
&&\quad\text{with solution $\bar{\mtx{M}}_{\star} := \bar{\mtx{B}}^\dagger \mtx{A}\bar{\mtx{B}}$;} \\
&\minimize_{\mtx{M}}\quad \fnorm{ \mtx{S}(\mtx{A}\bar{\mtx{B}} - \bar{\mtx{B}} \mtx{M}) }
&&\quad\text{with solution $\bar{\mtx{M}} := (\mtx{S}\bar{\mtx{B}})^\dagger (\mtx{SA}\bar{\mtx{B}})$.}
\end{aligned}
$$
Since $\mtx{Q}$ is an orthonormal basis for $\range(\bar{\mtx{B}})$,
we recognize that $\bar{\mtx{M}}_{\star}$ is similar to the RR matrix
$\mtx{Q}^*\mtx{AQ}$.  In view of~\cref{eqn:Mhat-form} and the relations
$\mtx{S}\bar{\mtx{B}} = \mtx{SBT}^{-1} = \mtx{U}$,
the sketched solution $\bar{\mtx{M}}$ is similar to the sRR matrix $\hat{\mtx{M}}$:
$$
\bar{\mtx{M}} = \mtx{U}^* (\mtx{SAB}) \mtx{T}^{-1}
	= \mtx{T} \hat{\mtx{M}} \mtx{T}^{-1}.
$$
Eigenvalue problems are invariant under similarity,
so it suffices to show $\bar{\mtx{M}} \approx \bar{\mtx{M}}_{\star}$.

To that end, we invoke~\cref{eqn:ls-residuals} columnwise to obtain the comparison
\begin{equation} \label{eqn:srr-anal-1}
\fnorm{ \mtx{A}\bar{\mtx{B}} - \bar{\mtx{B}} \bar{\mtx{M}}_{\star} }
	\leq \fnorm{ \mtx{A}\bar{\mtx{B}} - \bar{\mtx{B}} \bar{\mtx{M}} }
	\leq \frac{1+\eps}{1-\eps} \cdot \fnorm{ \mtx{A}\bar{\mtx{B}} - \bar{\mtx{B}} \bar{\mtx{M}}_{\star} }.
\end{equation}
Using the definitions of $\bar{\mtx{M}}_{\star}$ and $\mtx{Q}$, we find that
$$
\fnorm{ \mtx{A}\bar{\mtx{B}} - \bar{\mtx{B}} \bar{\mtx{M}}_{\star} }
	= \fnorm{ (\Id - \mtx{QQ}^*) \mtx{A} \bar{\mtx{B}} }
	\leq \sigma_{\max}(\bar{\mtx{B}}) \cdot \fnorm{ (\Id - \mtx{QQ}^*) \mtx{AQ} }.
$$
From the last two displays, a short argument using the triangle inequality
and the conditioning of the whitened basis produces
$$
\fnorm{ \bar{\mtx{M}} - \bar{\mtx{M}}_{\star} }
	\leq \frac{1}{\sigma_{\min}(\bar{\mtx{B}})} \cdot \fnorm{ (\bar{\mtx{M}} - \bar{\mtx{M}}_{\star} )\bar{\mtx{B}} }
	\leq \frac{2(1+\eps)}{(1-\eps)^2} \cdot \fnorm{ (\Id - \mtx{QQ}^*) \mtx{A} \mtx{Q} }.	 
$$
Here is an interpretation.  If $\range(\mtx{Q}) = \range(\mtx{B})$ is close
to an invariant subspace of $\mtx{A}$, then $(\Id - \mtx{QQ}^*) \mtx{A} \mtx{Q} \approx \mtx{0}$.
In this case, $\bar{\mtx{M}} \approx \bar{\mtx{M}}_{\star}$.
Therefore, sRR and RR solve nearby eigenvalue problems,
and we deduce that sRR is a backward stable approximation
to RR in exact arithmetic.

More generally, as long as $\range(\mtx{B})$ contains an approximate eigenvector of $\mtx{A}$ with small residual, \cref{eq:rescompare} shows that the same vector yields a comparably small residual for the sketched eigenproblem~\cref{eqn:srr-eig}.  Unfortunately, even in this case, there is no guarantee that sRR will find an approximate eigenpair with a small residual.  Indeed, the behavior of classic RR is already complicated, with pathological examples~\cite[p.~282]{stewart2}.  Nevertheless, RR is known to provide excellent outputs in the vast majority of cases; see~\cite{nakatsukasa2020sharp,Saa11:Numerical-Methods,Stewart99ageneralization} for the analysis.

\subsection{sRR with a Krylov subspace basis}

When $\mtx{B}$ is a (graded) basis for a Krylov subspace,
the analysis of sRR simplifies further.
In this case, the solutions $\mtx{M}_{\star}$ and $\hat{\mtx{M}}$
to~\cref{eqn:RR-matrix} and~\cref{eqn:sRR-matrix-form} differ only in the final column! 
To verify this point, observe that~\cref{eqn:RR-matrix} decouples into a family of $d$ least-squares problems:
$$
\minimize_{\mtx{m}_i \in \CC^{d}}\quad \fnorm{\mtx{A}\vct{b}_i - \mtx{B} \vct{m}_i}
\quad\text{for $i = 1, \dots, d$}.
$$
By construction, the vector $\mtx{Ab}_i$ lies in the span of $\mtx{B}$ for $i=1,\ldots, d-1$.
Each of these problems has a unique solution with zero residual.  Thus, the sketched problem~\cref{eqn:sRR-matrix-form}
correctly identifies the \textit{exact} solution.

\subsection{The symmetric case}
\label{sec:symeig}

Consider the symmetric eigenvalue problem
\begin{equation} \label{eqn:symeig}
\mtx{A}\vct{x} = \lambda\vct{x}
\quad\text{where $\mtx{A} = \mtx{A}^*$.}
\end{equation}
We can apply sRR directly to~\cref{eqn:symeig}.
Unfortunately, sRR is not guaranteed to (and in fact does not always)
return real eigenvalue estimates. At root, the sketched eigenvalue problem~\cref{eqn:srr-eig}
is not (similar to) a symmetric problem.  Accordingly, the computed eigenvectors need not be orthogonal.
This is an inherent drawback.

Fortunately, for~\cref{eqn:symeig}, sRR often computes eigenvalue estimates
that are real (or nearly real), and the associated eigenvectors tend to be nearly orthogonal.
We can anticipate this outcome when the whitened matrices satisfy
$\bar{\mtx{M}} \approx \bar{\mtx{M}}_{\star}$. 
Indeed, the eigenvalues of a symmetric matrix are well-conditioned
under nonsymmetric perturbations~\cite{kahan1975spectra,stewart-sun:1990}.
As for the eigenvectors, if two approximate eigenpairs $(\hat{\vct{x}}_1,\hat\lambda_1)$ and $(\hat{\vct{x}}_2,\hat\lambda_2)$
of a symmetric matrix $\mtx{A}$ have small residuals and sufficient gap $|\hat\lambda_1-\hat\lambda_2|$,
then it follows that $\hat{\vct{x}}_1,\hat{\vct{x}}_2$ are nearly orthogonal~\cite{parlettsym}.
This forces the high-quality eigenvectors computed by sRR to be nearly orthogonal.

For a real symmetric matrix $\mtx{A}$, our implementation of sRR simply extracts the real part of the computed eigenvalues and eigenvectors. 
When $\mtx{A}$ is complex Hermitian, we force the eigenvalues to be real (but not eigenvectors).
The design of a fast algorithm that respects symmetry remains an open problem.

\section{Constructing a basis for sRR}
\label{sec:eigs-basis}

The performance of RR and sRR depends on the quality of
the basis construction.  For these problems, it is natural
to consider \textit{block} Krylov bases generated by random vectors.
As before, we can consider nonorthogonal basis constructions.
Owing to the overlap with the discussion of single-vector
Krylov spaces, our presentation here is more telegraphic.

\subsection{Block Krylov subspaces}\label{sec:blockkry}

For the eigenvalue problem~\cref{eqn:eig}, we can search for solutions
using sRR with a block Krylov subspace.  Let $\mtx{\Omega} \in \CC^{n \times b}$
be an initial matrix; the dimension $b$ is called the \textit{block size}.  Define 
$$
\set{K}_p(\mtx{A}; \mtx{\Omega})
	:= \lspan\{\mtx{\Omega}, \mtx{A\Omega}, \dots, \mtx{A}^{p-1}\mtx{\Omega} \}
	= \lspan\{ \varphi(\mtx{A}) \mtx{\Omega} : \deg(\varphi) \leq p - 1 \}.
$$
Setting $d = bp$, we can express a basis $\mtx{B} \in \CC^{n \times d}$
for this subspace in the form
$$
\mtx{B} = [\mtx{B}_1, \dots, \mtx{B}_p] = [ \varphi_1(\mtx{A})\mtx{\Omega}, \dots, \varphi_p(\mtx{A}) \mtx{\Omega} ].
$$
Here, $\{ \varphi_i : i = 1, \dots, d \}$ is a linearly independent family of filter polynomials.
For eigenvalue problems, the generating matrix $\mtx{\Omega} \in \CC^{n \times b}$
may be drawn at random from a standard normal distribution.\footnote{In this context, we do not derive much computational benefit
from fancier nonadaptive distributions, such as SRFTs or sparse embeddings.}

Historically, the NLA literature has prescribed a small block size, say $b \leq 4$,
and a large depth $p$.  More recent research~\cite[Sec.~11]{MartinssonTroppacta}
has identified an opportunity to use a large block size $b$, say 10s or 100s,
with a much smaller depth, say $p \leq 10$.  This shift in perspective has already
transformed the computational profile of block Krylov methods for low-rank approximation.
For instance, we can parallelize the computation
over the columns of $\mtx{\Omega}$ (or over the filter polynomials $\varphi_{i}$).
In combination with sRR, nonorthogonal basis constructions promise further benefits.
See~\cite{muscolanczos,Tropp_NM22,Tro18:Analysis-Randomized-TR} for theoretical analysis of block Krylov subspaces
for low-rank matrix approximation and symmetric eigenvalue problems.

\begin{remark}[Other kinds of bases]
There are other subspace projection methods for solving eigenvalue problems
that use bases other than Krylov subspaces.  For example, the Jacobi--Davidson method and the LOBPCG algorithm uses alternative ideas to build a search space. These methods may also be combined with sRR. 
\end{remark}

\subsection{Basis diagnostics and restarting}

As with single-vector Krylov subspaces, we can sketch basis vectors
as they are generated to collect summary information about the quality
of the basis.  Indeed, if $\mtx{B} \in \CC^{n \times d}$ is a basis,
then the condition number of the sketched basis $\kappa_2(\mtx{SB})$
serves as a proxy for $\kappa_2(\mtx{B})$; see~\cref{eqn:cond-diagnostic}.
When the basis is poor, it can be important to use stabilized sRR (\cref{sec:srr-stabilization}).

It can also be effective to restart production of the Krylov subspace
when the quality of the basis starts to decline.
For example, we may compute a basis for the Krylov subspace $\set{K}_p(\mtx{A}; \mtx{\Omega})$,
and we can feed this basis to sRR to extract a matrix $\mtx{X}$ whose columns
approximately span the desired invariant subspace of $\mtx{A}$.
Then we pass to the Krylov subspace $\set{K}_p(\mtx{A}; \mtx{X})$,
and so forth.
Randomized subspace iteration~\cite{halko2011finding}
is a simple version of this technique.

Another possibility for restarting is to deflate converged eigenpairs by working in their orthogonal complement. 
Convergence of Ritz pairs and loss of orthogonality are known to be tightly linked~\cite[Ch.~11]{parlettsym}.
Felicitously, sRR is able to identify such eigenpairs cheaply.  Optimizing the sRR restarting strategy is left as future work.

\subsection{Block monomial basis with orthogonalization}

Although the monomial basis is anathema for large-degree polynomials,
we can still use it for shallow Krylov subspaces (say, when $p < 5$).
In this case, we can assemble a basis $\mtx{B} = [\mtx{B}_1, \dots, \mtx{B}_p]$
for $\set{K}_p(\mtx{A}; \mtx{\Omega})$ as follows.
Set $\mtx{B}_1 = \texttt{orth}(\mtx{\Omega})$, and iterate
$$
\mtx{B}_j = \texttt{orth}(\mtx{\Omega}_j)
\quad\text{where}\quad
\mtx{\Omega}_{j} = \mtx{A}\mtx{B}_{j-1}
\quad\text{for $j = 2, 3, \dots, p$.}
$$
We acquire the reduced matrix $\mtx{AB}$ as a by-product of this computation.

The block monomial basis has been used in the ``blanczos''
method~\cite{rokhlin2009randomized,HMST11:Algorithm-Principal,muscolanczos,Tro18:Analysis-Randomized-TR,MartinssonTroppacta}
for low-rank matrix approximation,
but it requires an expensive full orthogonalization of
$\mtx{B}$ in the final step.
When used as an input to sRR, it may not be necessary to
reorthogonalize the block monomial basis $\mtx{B}$.

\subsection{Block Arnoldi with truncation}

We can mitigate the rapid condition number growth of the block monomial
basis by adding extra orthogonalization steps.  For recurrence
length $k \in \mathbb{N}$, we set $\mtx{B}_1 = \texttt{orth}(\mtx{\Omega})$ and iterate
$$
\mtx{B}_j = \texttt{orth}(\mtx{\Omega}_j)
\quad\text{where}\quad
\mtx{\Omega}_j = (\Id - \mtx{B}_{j-1} \mtx{B}_{j-1}^* - \dots - \mtx{B}_{j-k} \mtx{B}_{j-k}^*) (\mtx{AB}_{j-1}).
$$
The resulting basis $\mtx{B} = [\mtx{B}_1, \dots, \mtx{B}_p]$ serves as an input to sRR.
The choice $k = 1$ or $k = 2$ already improves substantially over the block monomial basis. 

When $\mtx{A}$ is Hermitian, the choice $k = 2$ corresponds to
the block Lanczos method without reorthogonalization~\cite[Chap.~13]{parlettsym}.
Historically, the reorthogonalization step has been regarded as important
for achieving robustness.  If we use block Lanczos with sRR,
then we can often dispense with reorthogonalization.

As with sGMRES, in the unblocked case ($b=1$), we recommend $k$-truncated Arnoldi with a modest $k$ as shown in~\cref{alg:srr}. However, for eigenvalue computations, there is compelling reason to take the block size $b\gg 1$.  As the block size $b$ increases, the cost of orthogonalization quickly becomes devastating.

\subsection{Block Chebyshev recurrence}\label{sec:blockcheb}

By employing other polynomial recurrences, we can potentially
\textit{eliminate} all expensive computations that involve high-dimensional basis vectors.  In particular,
the shifted-and-scaled Chebyshev recurrence emerges as an appealing option.

Suppose that we have prior knowledge that the spectrum of $\mtx{A}$ is contained
in the axis-aligned rectangle $[c \pm \delta_x, \pm\delta_y]$, and set $\varrho = \max\{\delta_x,\delta_y\}$.
Then we can form a block Chebyshev basis $\mtx{B} = [\mtx{B}_1, \dots, \mtx{B}_p]$
for $\set{K}_p(\mtx{A}; \mtx{\Omega})$ as follows.  $$
\begin{aligned}
\mtx{B}_1 = \mtx{\Omega}; \qquad \mtx{B}_2 = \frac{1}{2\varrho} (\mtx{A} - c\Id) \mtx{\Omega}; \qquad
\mtx{B}_j = \frac{1}{\varrho} \left[ (\mtx{A} - c\Id) \mtx{B}_{j-1} - \frac{\delta_x^2 - \delta_y^2}{4\varrho} \mtx{B}_{j-2} \right].
\end{aligned}
$$
To implement this procedure, we typically need to perform a coarse initial eigenvalue computation
(using sRR + block Arnoldi) to obtain a rough estimate for the spectrum of $\mtx{A}$.
For this purpose, a small block size $b$ and depth $p$ usually suffice.
We may also consider Chebyshev polynomials based on rotated ellipses, as in~\cite{philippe2012generation}.

A remarkable feature of this approach is that we can compute the block Chebyshev
basis for $\set{K}_p(\mtx{A}; \mtx{\Omega})$ with $b(p-1)$ matvecs plus $O(nbp)$ operations. 
In contrast, it requires $O(n(bp)^2)$ extra operations to produce
an (approximately) orthogonal basis.
Beyond that, the Chebyshev recurrence can be implemented efficiently
in parallel or with SIMD processors, and the lack of inner products and
orthogonalization steps allows us to evade communication and synchronization costs.

\section{Computational experiments}
\label{sec:experiments}

This section presents numerics
that showcase the potential of sGMRES and sRR for
solving large linear systems and eigenvalue problems.
All examples involve real-valued  matrices, with appropriate
modifications to the methodology.
All computations were performed in MATLAB version 2020a on a workstation
with 256GB memory and 96 cores, each clocked at 3.3 GHz.

\subsection{Solving linear systems with sGMRES}

This subsection applies the sGMRES method to solve symmetric
and nonsymmetric linear systems.

\subsubsection{Algorithm details}

Our implementation of sGMRES follows the pseudocode in
\cref{alg:sgmres}. We construct a basis using $k$-truncated Arnoldi
with small values $k\in\{2,4\}$, unless otherwise noted.
In one example, we consider a Chebyshev basis, as described in
\cref{sec:chebyshev}.
We do not whiten the basis or restart sGMRES.
The subspace embedding is based on an SRFT matrix~\cref{eqn:SRFT}
where $\mtx{E}$ has independent Rademacher\footnote{A \textit{Rademacher random variable} takes values $\pm 1$ with equal probability.}
entries and $\mtx{F}$ is a discrete cosine transform (DCT2).
This sketch is easy to implement, but it uses $O(nd\log n)$
operations rather than $O(nd \log d)$.

We do not report tests involving the more elaborate algorithms
discussed in \cref{sec:sgmres-algs}, because fine-tuning for optimal performance
is outside the scope of this exploratory research.

\subsubsection{A nonsymmetric linear system}\label{sec:nonsymlin}

This subsection offers details about solving the nonsymmetric linear system,
documented in Figure~\ref{fig:gmres-intro} of the introduction.
The matrix $\mtx{A}$ is the sparse instance \texttt{t2em}
with dimension $n = 921,632$ from the SuiteSparse Matrix Collection~\cite{davis2011university}.
The right-hand side $\vct{f}$ is generated via $\vct{f}=\mtx{A}\vct{x}$, where $\vct{x}$ is a random vector drawn from the standard normal distribution.
We compare sGMRES with the MATLAB command
\texttt{gmres} without restarting and with restarting frequencies $\{10,30,100\}$.  Observe that
the $k$-truncated Arnoldi basis leads to a reduced matrix $\mtx{AB}$ whose condition number
grows quickly, but the condition number remains below the tolerance $u^{-1}$ throughout the computation.
This property ensures that sGMRES is effective.

\begin{figure}[t] \includegraphics[width=0.50\textwidth]{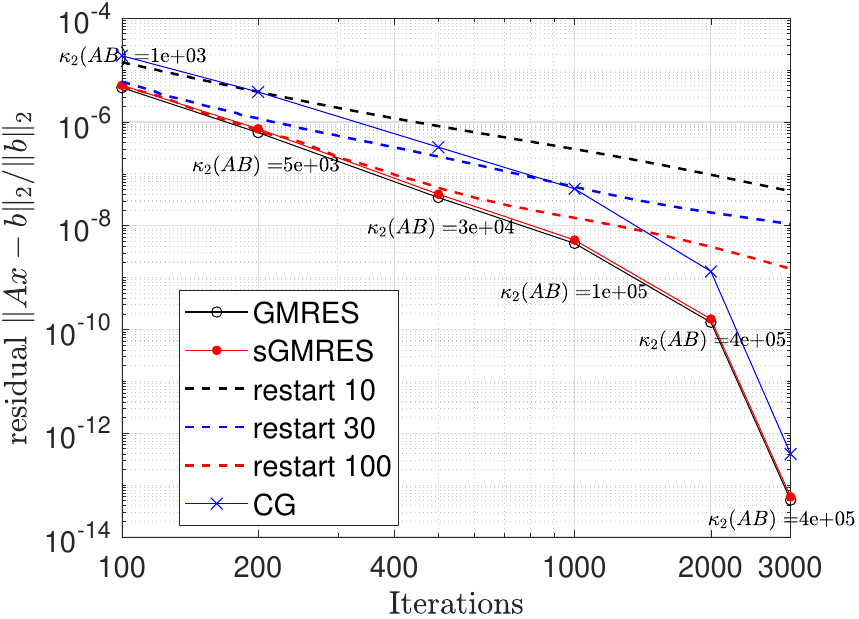}  
\includegraphics[width=0.47\textwidth]{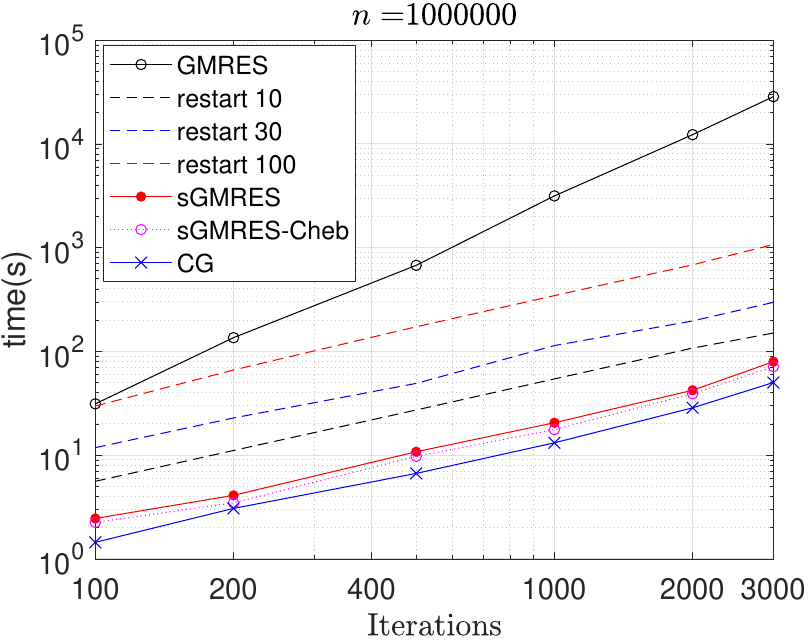}  
\caption{
\textbf{sGMRES versus GMRES: Laplacian system.}  These panels compare the performance of MATLAB \texttt{gmres}
(with and without restarting) against the sGMRES algorithm (with $2$-partial orthogonalization or the Chebyshev basis).
The sparse linear system $\mtx{A}\vct{x} = \vct{f}$ involves a 2D Laplacian matrix with dimension $n = 10^6$.
\textbf{Left:} Relative residual and condition number $\kappa_2(\mtx{AB})$ of the reduced matrix associated with
the $k$-truncated Arnoldi basis.  \textbf{Right:} Total runtime including basis generation.
} \label{fig:gmres-lap}
\end{figure}

\subsubsection{A symmetric linear system}\label{sec:linsyslap}

We consider a symmetric test matrix $\mtx{A}$ with dimension $n = 10^6$,
obtained by discretizing the 2D Laplacian.\footnote{To generate the matrix we used the code in
\url{https://www.mathworks.com/matlabcentral/fileexchange/27279-laplacian-in-1d-2d-or-3d}.}
The matrix is positive semidefinite with kernel $\mathbf{e} = [1, 1, \dots, 1]^*$.
We solve the Poisson problem $\mtx{A}\vct{x} = \vct{f}$ where
the right-hand side $\vct{f}=\mtx{A}\vct{x}$ is generated as above, which forces $\vct{e}^*\vct{f}=0$. 
For this problem, CG would be more appropriate than GMRES.
In fact, specialized algorithms (e.g., multigrid) are available,
but this example still offers an inspiring illustration of our methodology.

Figure~\ref{fig:gmres-lap} describes the progress of sGMRES where the basis is generated
by $k$-truncated orthogonalization for $k = 2$ and where the basis is generated by the
Chebyshev recurrence (more below).  We compare with the MATLAB commands \texttt{pcg}
and \texttt{gmres} without restart and with restart frequencies in $\{10,30,100\}$.
For all methods, the initial solution $\vct{x}_0 = \vct{0}$.

For both versions of sGMRES, the cost of $d$ iterations is 
about 50$\%$ slower than $d$ iterations of CG.  Nevertheless, for the same number $d$ of iterations,
sGMRES achieves $\ell_2$ residual norms that are about $5\times$ smaller than CG.
According to this metric, the sGMRES method is more efficient than CG.  As we saw for the
nonsymmetric problem, sGMRES is up to $6\times$ faster than the restarted versions
of GMRES, which do not converge to high accuracy.  Meanwhile, sGMRES achieves
the same accuracy as GMRES, but the sketched version is up to $100\times$ faster
after $3000$ iterations.

The Laplacian matrix is a natural candidate for testing the Chebyshev basis
because we have prior knowledge about the spectrum.  We use the fact that
the eigenvalues are real numbers in the interval $[0,8]$ to select the
parameters for the Chebyshev recurrence (\cref{sec:chebyshev}).  The Chebyshev basis construction
is slightly faster than the $k$-truncated Arnoldi construction because it
requires no inner products or orthogonalization steps.  Even so, the
quality of the Chebyshev basis is decent; after $3000$ iterations,
the reduced matrix has condition number $\kappa_2(\mtx{AB}) \approx 10^8$,
which is good enough for sGMRES to succeed.  
The $k$-truncated Arnoldi basis is still better conditioned (see Figure~\ref{fig:symeig}), but we do not need this improvement.
This experiment is intriguing because the Chebyshev basis can offer dramatic
benefits in parallel computing environments~\cite{JC91:Parallelizable-Restarted,philippe2012generation,ballard2014communication,chen2020predict}.

\subsection{sGMRES: Hard examples}

It is important to acknowledge that the sGMRES method is not
always an effective tool for solving linear systems.
In some cases, sGMRES inherits its weaknesses from
GMRES, but there are also new phenomena that arise.

First, there are linear systems where classic GMRES
cannot produce a small residual because the Krylov subspace
does not have sufficient approximation power.  sGMRES cannot
cure this debility.  In these cases, preconditioning
is critical.

Second, sGMRES is not especially useful for problems where the matrix--vector multiply $\vct{x} \mapsto \mtx{Ax}$
is costly relative to the other arithmetic.
For example, when $\mtx{A}$ is dense, over 99\% of the
runtime of GMRES or sGMRES may be devoted to matvecs.

Third, and most seriously, there are linear systems where
it is very difficult to construct a numerically full-rank basis
for the Krylov subspace without meticulous orthogonalization.
The rest of this subsection documents one such problem
instance. 

\begin{figure}[t] \includegraphics[width=0.47\textwidth]{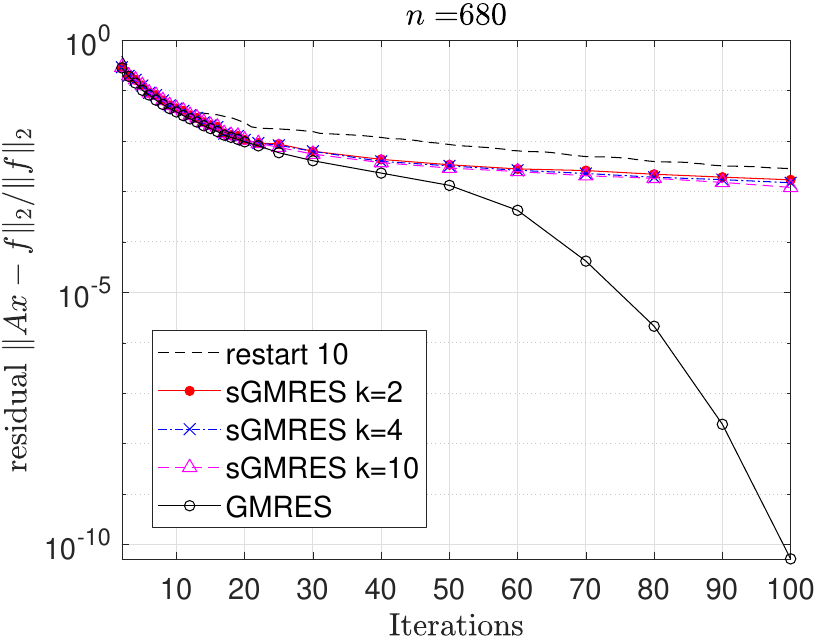}  
\includegraphics[width=0.50\textwidth]{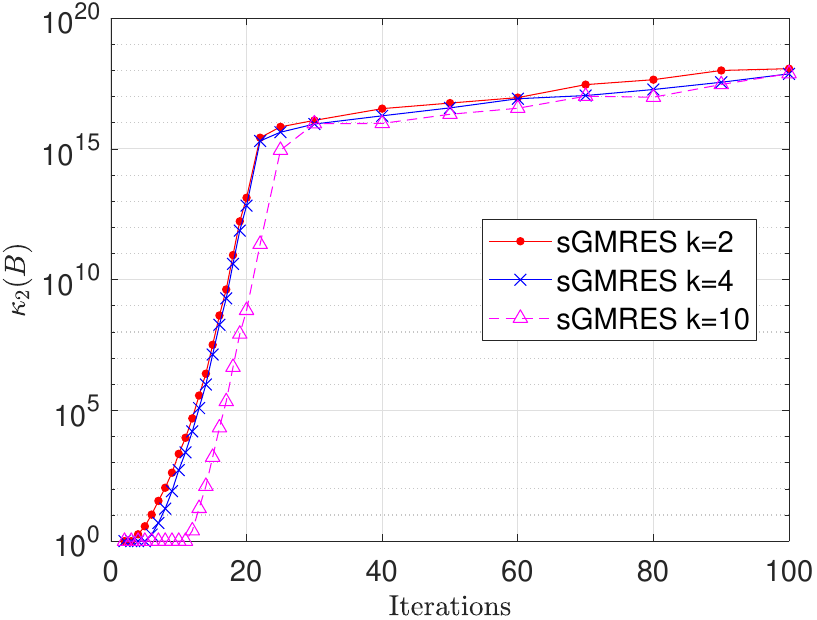}  
\caption{\textbf{sGMRES: Hard problems.}  For some linear systems, it is expensive to construct a well-conditioned
basis $\mtx{B}$ for the Krylov subspace.  When $\kappa_2(\mtx{AB}) > u^{-1}$, the sGMRES algorithm may fail to match
the GMRES algorithm with full orthogonalization.  \textbf{Left:} Relative residual norms. \textbf{Right:} Condition number $\kappa_2(\mtx{B})$ of the $k$-truncated Arnoldi basis.  When $\kappa_2(\mtx{B})>u^{-1}$, the reported values are unreliable. 
}\label{fig:gmres-fail}
\end{figure}

Consider the matrix \texttt{FS 680 1} from Matrix Market,
which is known to instigate Krylov bases with bad behavior~\cite[Table 2]{philippe2012generation}.
In this case, the basis $\mtx{B}$ and the reduced matrix $\mtx{AB}$
and their sketches $\mtx{SB}$ and $\mtx{SAB}$ have rapidly
increasing condition number.
Once $\kappa_2(\mtx{SAB}) > u^{-1}$, numerical errors can cause sGMRES to fail,
even when GMRES is successful.
See Figure~\ref{fig:gmres-fail} for an illustration,
which shows that increasing the extent $k$ of the truncation does not help.

We can always monitor the conditioning of the reduced matrix $\mtx{AB}$
inexpensively by means of its sketch $\mtx{SAB}$.  Unfortunately, we are not aware of a reliable
mechanism for controlling the conditioning, short of full
orthogonalization.  Indeed, $k$-truncated Arnoldi does not
even guarantee monotone decrease of the condition number as $k$ increases.
This issue remains a challenge for sGMRES.
The ideas from~\cite{balabanov2020randomized} may be useful here.

\subsection{Solving eigenvalue problems with sRR}

This subsection studies the performance of sRR for solving
nonsymmetric and symmetric eigenvalue problems.

\subsubsection{Algorithm details} \label{sec:srr-details}

Our implementation of sRR follows the pseudocode in~\cref{alg:srr}
with minor changes to facilitate comparison with the MATLAB \texttt{eigs} command. In particular, we focus on a single-vector Krylov subspace ($b = 1$),
and we use $k$-truncated Arnoldi to form the basis.
In one example, we consider a block Krylov subspace with a Chebyshev
basis, as described in \cref{sec:blockcheb}.
We do not use restarting or stabilization, except as noted.
The subspace embedding is based on an SRFT matrix~\cref{eqn:SRFT}
where $\mtx{E}$ is diagonal Rademacher and $\mtx{F}$ is a DCT2.

When using \texttt{eigs}, we set the option
\texttt{opts.p=r; opts.maxit=1;}
to suppress restart.  We set the flag \texttt{opts.issym}
to reflect whether the problem is symmetric.

For eigenvalue computations, we anticipate that block Krylov subspaces
can yield significant advantages over single-vector Krylov subspaces.
Some of these improvements derive from higher-order BLAS.  We can
also take advantage of SIMD architectures, and we can reduce costs of
communication and synchronization in parallel computing environments.

\subsubsection{Nonsymmetric eigenvalue problems}\label{sec:nonsymeig}

This section describes the nonsymmetric eigenvalue problem
that forms the basis for Figure~\ref{fig:eigs-intro}.  This
computation is modeled on the trust-region subproblem (TRS) from optimization~\cite{conn2000trust}:
\begin{equation} \label{eqn:trs}
\minimize_{\mtx{x} \in \RR^n}\quad  \frac{1}{2}\mtx{x}^* \mtx{Ax}+\mtx{g}^*\mtx{x}\quad
\subjectto\quad \norm{\vct{x}}_2 \leq \Delta.
\end{equation}
This quadratic program can be reduced to a nonsymmetric eigenvalue problem~\cite{adachi2017solving}:
\begin{equation}  \label{eq:trsgep}
  \begin{bmatrix}
 \mtx{A} & \Delta^{-2} \vct{gg}^* \\-\Id & \mtx{A} 
  \end{bmatrix}\mtx{x} =   \lambda \mtx{x}.
\end{equation}
To obtain a solution to~\cref{eqn:trs}, we extract a (scaled) eigenvector of~\cref{eq:trsgep}
corresponding to the right-most eigenvalue (which must be real).

We consider an instance of~\cref{eq:trsgep} where $\mtx{A}$ is an $n \times n$ tridiagonal
matrix with equispaced values in $[-1,1]$ on the main diagonal and with $1$s on the off-diagonals.
In this case, the block matrix admits a fast matrix--vector multiplication operation.
The vector $\vct{g} \in \RR^n$ is drawn from the standard normal distribution and scaled
so that $\norm{\vct{g}}_2 = 0.01$.  The constraint value $\Delta = 1$.

We solved~\cref{eq:trsgep} using sRR, as described in \cref{sec:srr-details}, taking $k=2$.
The Krylov subspace was generated from the initial vector $\vct{0} \oplus \vct{g}$.
The results appear in~\cref{fig:eigs-intro}.  The modest loss of accuracy
in sRR after $1500$ iterations can be remedied by using the stabilization
process (\cref{sec:srr-stabilization}), which approximately doubles
the runtime.

\begin{figure}[t] \includegraphics[width=0.46\textwidth]{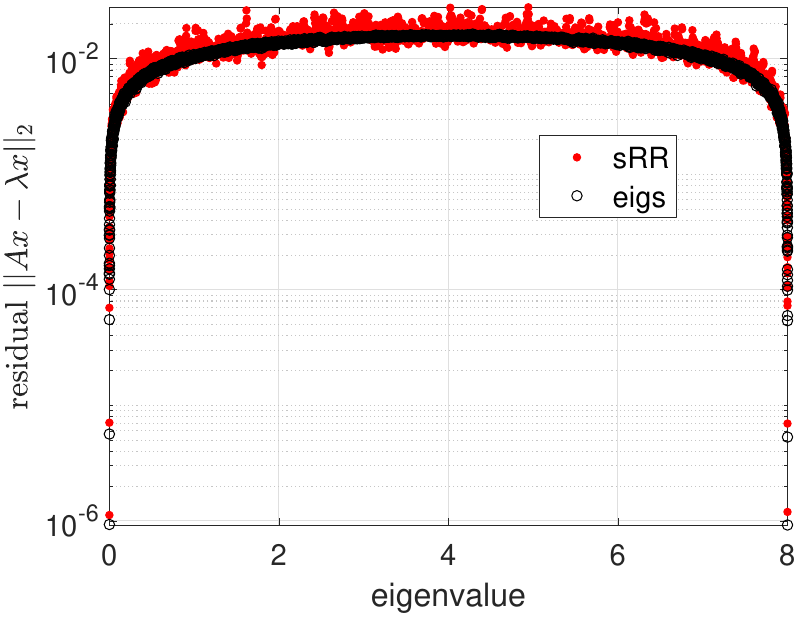}  
\includegraphics[width=0.48\textwidth]{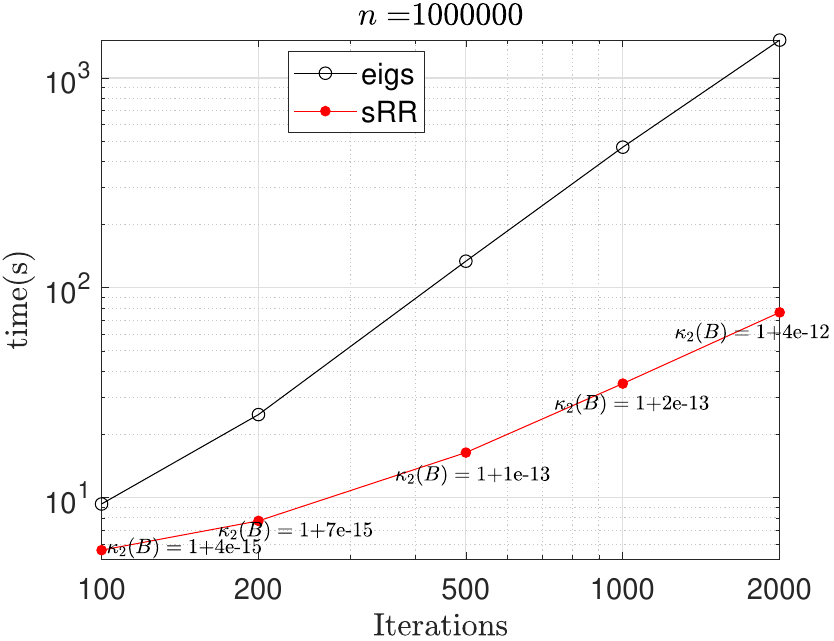}  
\caption{\textbf{sRR versus RR: Symmetric eigenvalue problem.}
These panels compare the performance of MATLAB \texttt{eigs} (without restarting) against the sRR algorithm (where the basis $\mtx{B}$
is computed by the Lanczos recurrence).  The sparse, symmetric eigenvalue problem
$\mtx{A}\vct{x} = \lambda \vct{x}$ has dimension $n = 10^6$, and it arises from a 2D Laplacian.
\textbf{Left:} Residuals as a function of computed eigenvalues after 2000 iterations.
\textbf{Right:} Total runtime, including basis generation, and the condition number $\kappa_2(\mtx{B})$ of the basis.
} \label{fig:symeig}
\end{figure}

\subsubsection{Symmetric eigenvalue problems}

Next, we present an example of a symmetric eigenvalue problem.  The matrix $\mtx{A}$ is the 2D Laplacian matrix with dimension $n = 10^6$ that was introduced in~\cref{sec:linsyslap}.  The initial vector for the Krylov subspace is drawn from the standard normal distribution, and it is shared between sRR and \texttt{eigs}.  We use the basic Lanczos recurrence  (i.e., $2$-truncation without extra orthogonalization) to construct the subspace basis.  For sRR, we report the real parts of the eigenvalues and eigenvectors, as discussed in~\cref{sec:symeig}.

\Cref{fig:symeig} displays the results of the experiment.  We see that sRR identifies the same eigenvalues as RR, and the sRR residual norms are within a small factor of the RR residual norms.  Completing 2000 iterations, sRR runs $12\times$ faster than \texttt{eigs}.  The difference is likely because \texttt{eigs} enforces orthogonality to ensure that the Lanczos method remains robust.  In this instance, the basic Lanczos method already constructs a basis that is almost orthogonal, but sRR does not require this property to succeed.

When the Lanczos method is used to reduce the matrix to (partial) tridiagonal form, it is critical that the Lanczos basis remain almost perfectly orthogonal.  Loss of orthogonality of the basis leads to \textit{ghost eigenvalues}, which are repeated estimates of a single eigenvalue~\cite[Ch.~7]{demmelbook}.  (Selective) orthogonalization is a traditional remedy~\cite{simon1984analysis,simon1984lanczos}, but it can be costly.  In our experience, sRR rarely produces ghost eigenvalues because it does not need the basis to reduce the matrix to tridiagonal form.

\begin{figure}[t] \centering
  \begin{minipage}[t]{0.49\hsize}
\includegraphics[width=.95\textwidth]{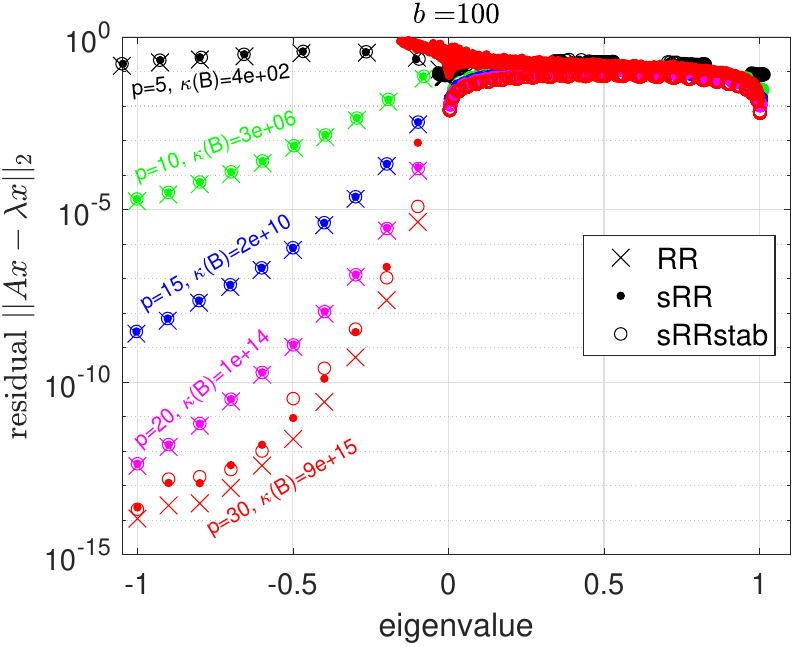}
  \end{minipage}   
\begin{minipage}[t]{0.49\hsize}
\includegraphics[width=.95\textwidth]{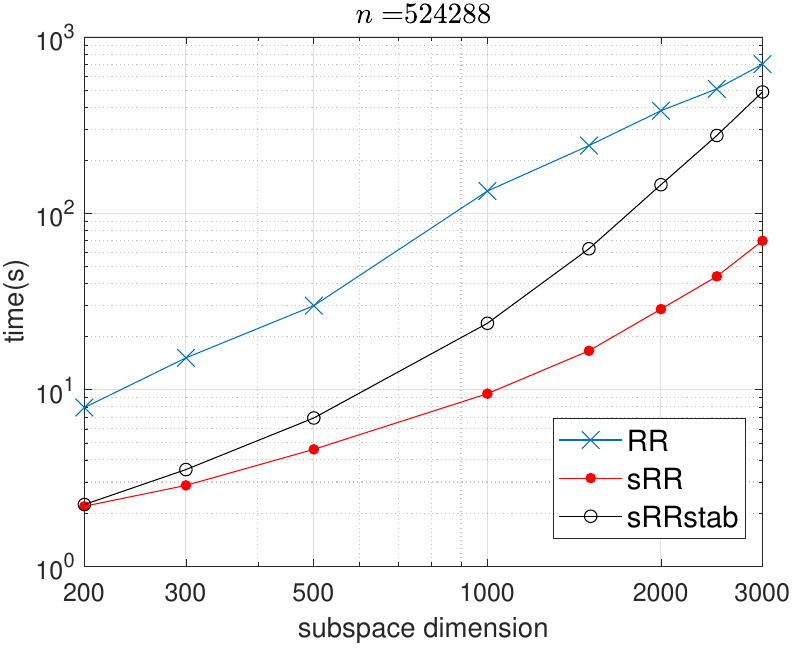}
  \end{minipage}
\includegraphics[width=.8\textwidth]{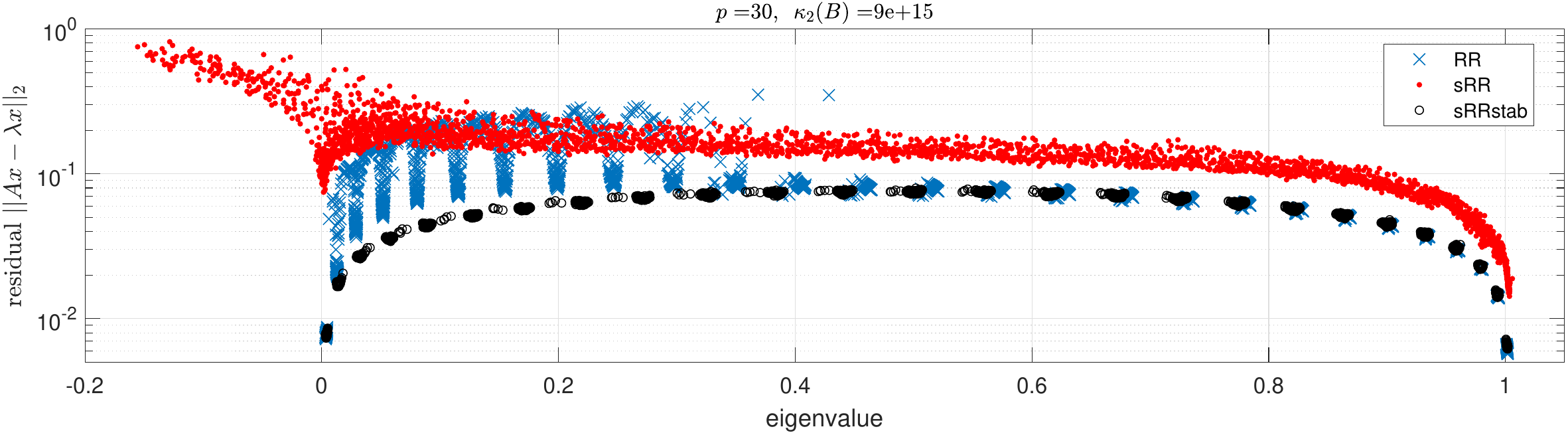}
  \caption{\textbf{sRR: Ill-conditioned bases and stabilization.}  This diagram shows how eigenvalue residuals improve
  with the depth $p$ of a block Krylov--Chebyshev subspace $\mtx{B} \in \RR^{n \times (bp)}$ with block
  size $b = 100$.
  \textbf{Left:} Residuals as a function of eigenvalue estimates, along with the condition number $\kappa_2(\mtx{B})$
  of the basis.
  \textbf{Right:} Runtime for eigenvalue extraction, excluding basis generation.
   \textbf{Bottom:} Magnification of the left panel for $p=30$ to illustrate (interior) eigenvalue estimates in $[0,1]$.
}
  \label{fig:rankdef}
\end{figure}

\subsubsection{Poorly conditioned bases and stabilization} \label{sec:expillcond}

When the computed basis $\mtx{B}$ is ill-conditioned, sRR may not produce reliable
eigenvalues estimates.  We can partially stanch this loss by stabilization, as discussed
in~\cref{sec:srr-stabilization}.

Without loss of generality, we consider a diagonal matrix $\mtx{A}$.  The dimension $n = 2^{19}$,
and the eigenvalues are ten equispaced spaced numbers in $[-1, -0.1]$, along with $2^{19} - 10$
equispaced numbers in $[0, 1]$.  Via the block Chebyshev recurrence (\cref{sec:blockcheb}),
we construct a (nonorthogonal) basis
$\mtx{B} \in \RR^{n \times (bp)}$ for the block Krylov subspace
with block size $b = 100$ and increasing depth $p$.
We then use classic RR, sRR, and sRRstab (the stabilized version) to compute
eigenpairs of $\mtx{A}$.  

\Cref{fig:rankdef} displays the results.  The condition number $\kappa_2(\mtx{B})$ of the basis grows
with the depth $p$ of the block Krylov subspace.  Even so, sRR computes
accurate eigenpairs while $\kappa_2(\mtx{B}) \lesssim u^{-1}$.  When the condition
number of the basis is larger, sRRstab even produces more accurate estimates for the
interior eigenpairs than RR does.
Excluding the cost of basis generation, sRR runs up to $15\times$ faster than classic RR.
The stabilized algorithm usually runs faster than RR, but the relatively large constant in the $O(d^3)$ cost of the
SVD and QZ algorithms dominates the runtime as the basis dimension $d$ grows.

\section{Variations and extensions}\label{sec:extension}

The ideas underlying sGMRES and sRR can be adapted to address a wide variety
of eigenvalue and singular value computations.

\subsection{Generalized eigenvalue problems}

Consider the problem
\begin{equation} \label{eqn:gep}
\text{Find nonzero $\vct{x} \in \CC^n$ and $\lambda \in \CC$ : }\quad
\mtx{H} \vct{x} = \lambda \mtx{J} \vct{x}
\quad\text{where $\mtx{H}, \mtx{J} \in \CC^{n \times n}$.} \end{equation}
Suppose that $\mtx{B} \in \CC^{n \times d}$ is a basis
that captures approximate solutions to~\cref{eqn:gep}.
Following the development in~\cref{sec:rr-faces}, the
classic RR method can be interpreted
as a variational problem:
\begin{equation} \label{eqn:gep-var}
\minimize_{\mtx{M} \in \CC^{d \times d}}\quad
\fnorm{ \mtx{HB} - \mtx{JBM} }.
\end{equation}
Given a solution $\mtx{M}_{\star} = (\mtx{JB})^\dagger(\mtx{HB})$
to~\cref{eqn:gep-var}, we pose the ordinary eigenvalue problem
$\mtx{M}_{\star}\vct{y} = \theta \vct{y}$.  Each eigenpair
$(\vct{y}, \theta)$ induces an approximate
solution $(\mtx{B}\vct{y}, \theta)$ to~\cref{eqn:gep}.

Given a subspace embedding $\mtx{S} \in \CC^{s \times n}$
for $\range( [\mtx{HB}, \mtx{JB}] )$, we can pass to the sketched problem
\begin{equation} \label{eqn:sgep-var}
\minimize_{\mtx{M} \in \CC^{d \times d}}\quad
\fnorm{ \mtx{S}(\mtx{HB} - \mtx{JBM}) }.
\end{equation}
The solution $\hat{\mtx{M}} = (\mtx{SJB})^\dagger (\mtx{SHB})$.
Then frame the ordinary eigenvalue problem $\hat{\mtx{M}} \vct{y} = \theta \vct{y}$.
Each eigenpair $(\hat{\vct{y}}, \hat{\theta})$ induces an approximate
solution $(\mtx{B}\hat{\vct{y}}, \hat{\theta})$ to~\cref{eqn:gep}.

Excluding basis generation, we can solve the generalized eigenvalue problem
via sketching with $O(d^3 + nd \log d)$ operations.  In contrast, the classic
RR approach typically requires $O(nd^2)$ operations.

\subsection{Low-rank matrix approximation} \label{sec:lowsvd}

The most successful application of randomized matrix computation
has been to approximate truncated singular value decompositions
efficiently~\cite{halko2011finding,MartinssonTroppacta}.
Using the new insights from our paper, we can accelerate these
algorithms by sketching.  The resulting techniques share
some genes with sketch-based algorithms for low-rank matrix
approximation~\cite{woodruff2014sketching,tropp2017practical,tropp2019streaming,nakatsukasa2020fast},
but they are different in spirit.

Let $\mtx{A} \in \CC^{m \times n}$ be a matrix.  Let $\mtx{B} \in \CC^{n \times d}$
be a basis, and suppose that we have access to the reduced matrix $\mtx{AB} \in \CC^{n \times d}$.
We can frame low-rank matrix approximation as a variational problem:
\begin{equation} \label{eq:lowrank-var}
\minimize_{\mtx{M}\in\CC^{d \times d}} \quad
\fnorm{ \mtx{ABM} - \mtx{A} }.
\end{equation}
The solution $\mtx{M}_{\star} = (\mtx{AB})^\dagger \mtx{A}$ produces the rank-$d$ matrix approximation
$$
\hat{\mtx{A}} = \mtx{AB} \mtx{M}_{\star} = (\mtx{AB})(\mtx{AB})^{\dagger} \mtx{A}
	= \mtx{QQ}^* \mtx{A},
$$
where $\mtx{Q} \in \CC^{n \times d}$ is an orthonormal basis for the range of $\mtx{AB}$.
If we choose $\mtx{B}$ at random, we obtain the Halko et al.~randomized SVD algorithm~\cite{halko2011finding}.
If we form an adapted basis $\mtx{B}$ by means of subspace iteration~\cite{rokhlin2009randomized,halko2011finding}
or block Krylov methods~\cite{rokhlin2009randomized,HMST11:Algorithm-Principal,muscolanczos,Tro18:Analysis-Randomized-TR},
we obtain much better approximations, as described in the cited work.

Let $\mtx{S} \in \CC^{s \times n}$ be an ``affine space'' embedding~\cite{Woo20:Algorithms-Big} with $s = 2d$.
(The SRFT~\cref{eqn:SRFT} and sparse map~\cref{eqn:sparse-map} both qualify.)
We pose the sketched problem
\begin{equation} \label{eqn:lowrank-sketch}
\minimize_{\mtx{M}\in\CC^{d \times d}} \quad
\fnorm{ \mtx{S}(\mtx{ABM} - \mtx{A}) }.
\end{equation}
The solution $\hat{\mtx{M}} = (\mtx{SAB})^\dagger (\mtx{SA})$
yields the rank-$d$ matrix approximation
\begin{equation} \label{eqn:lowrank-sketch-approx}
\hat{\mtx{A}}_{\mathrm{sketch}} = (\mtx{AB}) (\mtx{SAB})^\dagger (\mtx{SA}).
\end{equation}
The formula~\cref{eqn:lowrank-sketch-approx} is wholly unsuitable for practical computation,
but it can be replaced with a stable and efficient variant~\cite{nakatsukasa2020fast}.
If we choose $\mtx{B}$ to be a second sketching map, we obtain the (low-accuracy)
sketched SVD algorithms from~\cite{woodruff2014sketching,tropp2017practical,nakatsukasa2020fast}.

Our work delivers the novel insight that using an \textit{adapted} basis $\mtx{B}$
in~\cref{eqn:lowrank-sketch-approx} leads to a fast \textit{and} accurate algorithm for low-rank
matrix approximation.  Excluding the cost of basis generation,
we can stably form the approximation in $O(d^3 + (m+n) d\log d)$ operations.  In contrast with sketched
SVD algorithms, we attain errors similar to randomized subspace
iteration~\cite{halko2011finding} or randomized block Krylov
methods~\cite{muscolanczos,Tro18:Analysis-Randomized-TR}.

\section{Prospects}
\label{sec:prospects}

We believe that our framework for combining sketching
with subspace projection methods presents many exciting
opportunities and challenges.  Let us close by
highlighting some of the prospects.

First, our work suggests that traditional strategies for building
high-dimensional Krylov subspace bases merit a fresh look.  For
example, our experiments indicate that we can easily run thousands
of iterations of sGMRES, whereas orthogonalization dominates the
cost of classic GMRES after, say, a few dozen iterations.  One consequence
is that it would suffice to find a ``mediocre'' preconditioner
for linear systems that reduces the iteration complexity of sGMRES
to 1000s of iterations, rather than the historical goal of 10s of iterations.
Other aspects of basis generation that deserve further attention include
restarting, deflation, and pruning.

Second, we believe that the performance advantages of sGMRES and sRR
algorithms would be maximized in modern computing environments where
communication and synchronization costs dominate
computation~\cite{ballard2011minimizing}.  For example,
we can trivially parallelize the computation of block Krylov subspaces.
Likewise, sketching allows us to perform approximate orthogonalization of distributed
vectors by means of short messages.  
While our experiments focused on a serial computing environment,
there are clear opportunities for efficient implementations on GPUs,
multicore and parallel processors, distributed and cloud computing
systems, and so forth.

Third, aside from GMRES and RR, there are many
subspace projection methods that might benefit
from sketching.  For instance, there is an important class of algorithms (BiCG, BiCGstab, CGS, QMR, etc.)
for solving linear systems by means of Lanczos biorthogonalization.
These methods form Krylov subspaces with respect to both $\mtx{A}$
and $\mtx{A}^*$ using three-term recurrences, but they
have complicated stability properties.  Perhaps, with
sketching, we can improve the profile of these methods.

Finally, let us mention one remaining difficulty.  At present, we
lack a reliable mechanism for guaranteeing that the condition
number of basis $\mtx{B}$ and the reduced matrix $\mtx{AB}$
do not explode.  Truncated orthogonalization is a practical approach
that often works well, but it can fail.  It would
be valuable to identify strategies for inexpensively
producing computational bases that are numerically full rank.

\section*{Acknowledgments}

The authors would like to thank 
Oleg Balabanov, Alice Cortinovis,
Ethan Epperly, Laura Grigori, Ilse Ipsen, Daniel Kressner, Gunnar Martinsson, Maike Meier, Florian Schaefer, Daniel Szyld, Alex Townsend, Nick Trefethen, and Rob Webber
for valuable discussions and feedback.  We are grateful to Zden{\v e}k Strako{\v s} for sharing his insights
on Krylov subspace methods and the prospects for sGMRES and sRR.

\def\noopsort#1{}\def\l{\char32l}\def\v#1{{\accent20 #1}}
  \let\^^_=\v\def\hbk{hardback}\def\pbk{paperback}

\end{document}